# The Least Square-free Primitive Root Modulo a Prime

Morgan Hunter

October 2016

A thesis submitted in partial fulfilment of the requirements for the degree of
Bachelor of Science (Honours)
in Mathematics at the Australian National University.

Australian National University

# Declaration

The work in this thesis is my own except where otherwise stated.

Morgan Hunter

# Acknowledgements

I would really like to thank my supervisor Tim Trudgian. He suggested a great topic for my thesis that I have found very interesting and rewarding to study. I would also like to thank him for funding my MSI Honours Scholarship, without which I would not have been able to complete honours this year.

I would also like to thank Mark Lawrenson who assisted with the implementation of the algorithm.

# Abstract

The aim of this thesis is to lower the bound on square-free primitive roots modulo primes. Let $g^{\square}(p)$ be the least square-free primitive root modulo $p$. We have proven the following two theorems.

**Theorem 0.1.**

$$g^{\square}(p) < p^{0.88} \quad \textit{for all primes } p.$$

**Theorem 0.2.**

$$g^{\square}(p) < p^{0.63093} \quad \textit{for all primes } p < 2.5 \times 10^{15} \textit{ and } p > 9.63 \times 10^{65}.$$

Theorem 0.1 shows an improvement in the best known bound while Theorem 0.2 shows for which primes we can prove the theoretical lower bound.

After some introductory information in Chapter 1, we will start to prove the above theorems in Chapter 2. We will introduce an indicator function for primitive roots of primes in §2.1 and together with results from §1.2.1, §1.2.3 and §1.2.4 we will outline the first step in proving a general theorem of the above form. The next two stages in the proof will be outlined in Chapter 3. These two stages require the introduction of the prime sieve. Before defining the sieve in §3.2 we will introduce the $e-$free integers which will play an important role in defining the sieve.

In §3.2 we will obtain results that do not require computation, including Theorem 0.2. An algorithm is then introduced in §3.3 which is the last stage of the proof. There we will complete the proof of Theorem 0.1.

We are preparing to publish the results of this thesis.

# Contents

# 1 Introduction

There are many interesting questions concerning the distribution of primitive roots modulo primes. In particular we are interested in the least prime primitive root of a prime. The asymptotic bound on the least prime primitive root is quite weak, and very difficult to improve, so in this thesis we will instead concentrate on a more general type of primitive root. We will be studying the least square-free primitive root modulo a prime.

Before we are able to bound the least square-free primitive root, we need to understand what a primitive root is and what basic properties it has. After outlining this in §1.2 we will introduce some arithmetic functions that are important in later chapters. In §1.2.4 we introduce the famous Pólya–Vinogradov inequality. This inequality is crucial in lowering the bound on the least square-free primitive root of a prime.

## 1.1 Notation

Throughout this thesis standard analytic number theory symbols are used. We will use the following shorthand: $[x]$ denotes the integer part of $x$ and $(a, b)$ denotes the greatest common divisor of $a$ and $b$. We will use 'O'-notation and $\ll$ and $\gg$ symbols as follows: for functions $f(x)$ and $g(x)$ the notation $f(x) = O(g(x))$ and $f(x) \ll g(x)$ mean that there exists a positive constant $M$ such that $|f(x)/g(x)| < M$ when $x$ is sufficiently large. The relation $f(x) \gg g(x)$ is interpreted as $g(x) \ll f(x)$.

## 1.2 Primitive roots

In order to define what a primitive root is we need to define the order of an integer.

**Definition 1.1.** Let $(a, m) = 1$, then the **order** of $a \pmod{m}$ is the smallest integer $k$ such that $a^k \equiv 1 \pmod{m}$ and is denoted $\mathrm{ord}_m(a) = k$.

**Proposition 1.1** (Theorem 8.4 from [30]). *Let $k$ be a positive integer. If $ord_m a = d$ then*

$$ord_m\left(a^k\right) = \frac{d}{(d, k)}.$$

*In particular, there are $\phi(d)$ distinct powers of $a$ with order $d$.*

**Definition 1.2.** We say $r$ is a **primitive root** modulo $m$ (alternatively, $r$ is a primitive root of $m$) if $\mathrm{ord}_m(r) = \phi(m)$. Equivalently, $r$ is a primitive root modulo $m$ if it generates the set of integers coprime to $m$. In particular, for all $a$ such that $(a, m) = 1$ there exists $k$ such that $r^k \equiv a \pmod{m}$.

**Remark 1.1.** Proposition 1.1 implies that if $r$ is a primitive root modulo $m$ then $r^k$ is a primitive root modulo $m$ if and only if $(k, \phi(m)) = 1$. Therefore if we have found one primitive root modulo $m$, we can generate all other primitive roots of $m$.

**Example 1.1.** Let us find the primitive roots of 7. Since 7 is prime, order is defined for all positive integers less than 7. We have

$$2^1 \equiv 2,\ 2^2 \equiv 4,\ 2^3 \equiv 1$$

so 2 is not a primitive root of 7. Next we try 3,

$$3^1 \equiv 3,\ 3^2 \equiv 2,\ 3^3 \equiv 6,\ 3^4 \equiv 4,\ 3^5 \equiv 5,\ 3^6 \equiv 1$$

so 3 is a primitive root of 7. Now that we have a primitive root we can generate all other primitive roots, $\phi(7) = 6$ so $k = \{1, 5\}$. Therefore the only other primitive root modulo 7 is $3^5 \equiv 5$.

The fact that a primitive root modulo $m$ generates the set of all coprime integers to $m$ leads to the following definition.

**Definition 1.3.** Let $r$ be a primitive root modulo $m$. Then for all integers $a$ such that $(a, m) = 1$, we define the base $r$ **discrete logarithm** of $a$, to be the unique integer $k \in \{1, \ldots, \phi(m)\}$ such that

$$r^k \equiv a \pmod{m}.$$

We denote this $\text{ind}_r(a) = k$. The base $r$ discrete logarithm is also known as the base $r$ index.

Note that it follows directly from Remark 1.1 that if $g$ is a primitive root modulo $m$, then $a$ is a primitive root modulo $m$ if and only if $(\text{ind}_g(a), \phi(m)) = 1$. The discrete logarithm not only provides some useful notation, but also the discrete logarithm modulo $\phi(m)$ shares the basic properties of logarithms.

It may be tempting to assume that all integers have primitive roots, however this is not true.

**Example 1.2.**
Consider $m = 8$.
Here we have $\phi(8) = 4$ and from Definition 1.1 the order is only defined for coprime integers to 8. Therefore we are left to consider $r = 1, 3, 5$, and 7 as possible primitive roots. However

$$1^2 \equiv 3^2 \equiv 5^2 \equiv 7^2 \equiv 1 \pmod{8},$$

so $\text{ord}_8(r) = 2 \neq 4$ for all $r = 1, 3, 5, 7$ and therefore there are no primitive roots modulo 8.

The next question to ask is, what integers have primitive roots? We will not prove the following propostion about the existence of primitive roots, the proofs can be found in Chapter 8 of Rosen's Book [30] or Chapter 10 of Apostol's Book [1].

**Proposition 1.2.**

1. *There exist primitive roots for all primes.*
2. *Powers of* $2$*, except for* $1$*,* $2$ *and* $4$*, do not have primitive roots.*
3. *There exist primitive roots for all powers* $p^k$ *and* $2p^k$ *where* $p$ *is an odd prime and* $k \geq 1$*.*

As we can see, our example from above fits into category 2, as 8 is a power of 2 and therefore does not have any primitive roots. From this point onwards we will be focusing on primitive roots modulo primes.

By Definition 1.2 we know that $r$ is a primitive root modulo a prime $p$ if

$$\{r^k \pmod p) \mid k = 1, \dots, p-1\} = \{1, 2, \dots, p-1\}.$$

and we also have that every prime has exactly $\phi(p-1)$ primitive roots (Chapter 8.2 of [30]). As we mentioned at the start of this chapter, we are interested in the distribution of primitive roots modulo primes. However to study the primitive roots of an unknown prime $p$ we first need some background on arithmetic functions. These play an important part in later chapters of this thesis

### 1.2.1 The Möbius and Euler totient functions

As mentioned above, to study primitive roots we require some background information on particular arithmetic functions. Arithmetic functions are real or complex functions that are defined on the set of natural numbers. In this section we will look at the Möbius function and the Euler totient function. In §1.2.3 we will introduce Dirichlet characters, which are also arithmetic functions. The properties of Dirichlet characters will be important in Chapter 2.

**Definition 1.4.**
The **Möbius function** is defined by

$$\mu(n) = \begin{cases} 0 & \text{if } p^2 \mid n \text{ for any prime } p, \\ (-1)^k & \text{if } n \text{ is the product of } k \text{ distinct primes}, \\ 1 & \text{if } n = 1. \end{cases}$$

The **Euler totient function** is defined by

$$\phi(n) = \# \{k \in \mathbb{Z} \mid 1 \leq k \leq n, (k, n) = 1\}.$$

That is $\phi(n)$ is the number of integers less than $n$, coprime to $n$.

These two arithmetic functions will appear repeatedly in the following chapters and some of their important properties are stated below.

**Proposition 1.3.** *Both $\mu(n)$ and $\phi(n)$ are multiplicative.*
*That is if $(n, m) = 1$ then $\mu(nm) = \mu(n)\mu(m)$ and $\phi(nm) = \phi(n)\phi(m)$.*

*Proof.* See Chapter 6 in [30]. □

**Remark 1.2.** If $(a, b) > 1$ then $\mu(ab) = 0$. This follows because there exists $c$ such that $c \mid a$ and $c \mid b$ and therefore $c^2 \mid ab$.

**Proposition 1.4** (Sums over divisors)**.** *Let $n \geq 1$. Then*

1. $\sum_{d|n} \mu(d) = \begin{cases} 1 & \text{if } n = 1, \\ 0 & \text{if } n > 1. \end{cases}$

2. $\sum_{d|n} \phi(d) = n.$

*Proof.* See Theorem 264 and §16.2 in [15]. □

An example of Möbius inversion shows how these two arithmetic functions are related [1]. For $n \geq 1$ we have

$$\phi(n) = \sum_{d|n} \mu(d)\frac{n}{d}.$$

**Proposition 1.5.**

$$\phi(n) = n\prod_{p|n}\left(1 - \frac{1}{p}\right),$$

*where the product is over the distinct prime divisors of $n$.*

*Proof.* See Theorem 62 in [15]. □

An important application of the Möbius function is related to square-free integers.

### 1.2.2 Square-free integers

**Definition 1.5.** An integer $n$ is said to be **square-free** if it is the product of distinct prime factors.

For example, 42 is a square-free integer, $42 = 2 \cdot 3 \cdot 7$, while $56 = 2^3 \cdot 7$ is not a square-free integer as 2 is a repeated prime factor.
It is clear from the above definition that all primes are square-free hence square-free integers are a weak generalisation of the primes. Recall that $\mu(n)$ is 0 if $n$ is divisible by the square of a prime and $\pm 1$ otherwise. Hence one possible indicator function of square-free integers is

$$|\mu(n)| = \begin{cases} 1 & \text{if } n \text{ is square-free,} \\ 0 & \text{otherwise.} \end{cases}$$

Another characteristic equation for square-free integers is given by Shapiro [31]. First note that all integers $n$ can be expressed as $n = s^2 q$ where $s$ is an integer and $q$ is square-free. Therefore from Proposition 1.4 we have

$$\sum_{d^2|n} \mu(d) = \sum_{d|s} \mu(d) = \begin{cases} 1 & \text{if } s = 1, \\ 0 & \text{otherwise.} \end{cases}$$

If $s = 1$ then $n$ is square-free and so a characteristic equation for square-free integers is

$$\sum_{d^2|n} \mu(d) = \begin{cases} 1 & \text{if } n \text{ is square-free,} \\ 0 & \text{otherwise.} \end{cases} \tag{1.1}$$

Now consider the number of square-free integers less than or equal to $x$. We then have to consider the sum

$$\sum_{\substack{n \le x \\ n = \square - \text{free}}} 1 = \sum_{n \le x} |\mu(n)|$$

which is asymptotic to $6\pi^{-2}x$ (Theorem 8.2.1 in [31]). Not only do we have an implicit bound on the number of square-free integers less than $x$, Cipu gives the following explicit bounds.

**Lemma 1.1** (Lemma 4.2 in [3]). *If $x \geq 1$ then*

$$\sum_{n\leq x} |\mu(n)| = \frac{6}{\pi^2}x + P(x), \quad \textit{with}$$

*(a)* $-0.103229\sqrt{x} \leq P(x) \leq 0.679091\sqrt{x}$ *for* $x \geq 1$,

*(b)* $|P(x)| \leq 0.1333\sqrt{x}$ *for* $x \geq 1664$,

*(c)* $|P(x)| \leq 0.036438\sqrt{x}$ *for* $x \geq 82005$,

*(d)* $|P(x)| \leq 0.02767\sqrt{x}$ *for* $x \geq 438653$.

These explicit bounds on the number of square-free integers less than or equal to $x$ will be important in Chapter 2, where we obtain results on the least square-free primitive root modulo a prime. The next section introduces the Dirichlet characters which will also be important in Chapter 2.

### 1.2.3 Dirichlet characters

A Dirichlet character is a certain type of arithmetic function. They are important in the study of primitive roots, in particular they appear in the indicator function for primitive roots modulo a prime (2.1).

**Definition 1.6.** Let $q$ be a positive integer. Then a **Dirichlet character** modulo $q$ is a function $\chi : \mathbb{N} \to \mathbb{C}$ with the following properties:

1. $\chi$ is periodic modulo $q$, i.e. $\chi(n+q) = \chi(n)$ for all $n \in \mathbb{N}$.

2. $\chi$ is completely multiplicative, i.e. $\chi(nm) = \chi(n)\chi(m)$ for all $n, m \in \mathbb{N}$.

3. $\chi(n) \neq 0$ if and only if $(n, q) = 1$.

The character

$$\chi_0(n) = \begin{cases} 1 & \text{if } (n,q) = 1 \\ 0 & \text{if } (n,q) > 1 \end{cases}$$

is called the **principal character** modulo $q$.

**Proposition 1.6.** *Let $\chi$ be a Dirichlet character modulo $q$. Then the values of $\chi$ are either 0 or $\phi(q)^{th}$ roots of unity.*

*Proof.* From the definition, if $\chi(n) \neq 0$ then $(n, q) = 1$. By Euler's theorem (Theorem 5.15 in [30]) we have $n^{\phi(q)} \equiv 1 \pmod{q}$. Then as $\chi$ is multiplicative and periodic we have $\chi(n)^{\phi(q)} = \chi(n^{\phi(q)}) = \chi(1) = 1$. ∎

It follows from Proposition 1.6 that if $\chi$ is a Dirichlet character modulo $q$ then so is the complex conjugate $\overline{\chi}$, where $\overline{\chi}(n) = \overline{\chi(n)}$.

Just as we have defined the order of an integer modulo $q$ (Definition 1.1) we can define the order of a Dirichlet character modulo $q$. Let $\chi$ be a Dirichlet character modulo $q$ then the order of $\chi$ is the smallest exponent $d$, with $d \mid \phi(q)$, such that $\chi^d = \chi_0$.

**Proposition 1.7** (§6.5 in [31]). *There are exactly $\phi(q)$ Dirichlet characters modulo $q$. They are denoted $\chi_0$, $\chi_1$, $\ldots$, $\chi_{\phi(q)-1}$. In particular, given $d \mid \phi(q)$ there are $\phi(d)$ Dirichlet characters modulo $q$ of order $d$.*

Consider all the Dirichlet characters, $\chi$, modulo $q$. The possible orders of these characters are the divisors of $\phi(q)$. Let $d_1$, $d_2$, $\ldots$, $d_s$ be the divisors of $\phi(q)$. Then from Proposition 1.4 we have

$$\phi(d_1) + \phi(d_2) + \cdots + \phi(d_s) = \phi(q).$$

Since there are $\phi(q)$ Dirichlet character modulo $q$, Proposition 1.4 shows that the Dirichlet characters can be partitioned according to their order.

**Example 1.3.** When $q = 1$ or $q = 2$, $\phi(q) = 1$ and so there is only one Dirichlet character, namely the principal character.
When $q = 3$ or $q = 4$ then there are 2 Dirichlet characters defined in Table 1.1 and Table 1.2.

| $n$ | 1 | 2 | 3 |
|---|---|---|---|
| $\chi_0(n)$ | 1 | 1 | 0 |
| $\chi_1(n)$ | 1 | $-1$ | 0 |

Table 1.1: $q = 3$, $\phi(3) = 2$

| $n$ | 1 | 2 | 3 | 4 |
|---|---|---|---|---|
| $\chi_0(n)$ | 1 | 0 | 1 | 0 |
| $\chi_1(n)$ | 1 | 0 | $-1$ | 0 |

Table 1.2: $q = 4$, $\phi(4) = 2$

When $q = 5$ we have 4 Dirichlet characters. When $(n, 5) = 1$ Proposition 1.6 shows that the possible values of $\chi(n)$ are $\pm 1$ and $\pm i$. Since $\chi$ is multiplicative we have $\chi(2)\chi(3) = \chi(6) = \chi(1) = 1$. Also $\chi(4) = \chi(2)^2$ and so we can define the Dirichlet characters modulo 5 in Table 1.3. The last example is when $q = 6$. Once again there are only 2 Dirichlet characters, defined in Table 1.4.

| $n$ | 1 | 2 | 3 | 4 | 5 |
|---|---|---|---|---|---|
| $\chi_0(n)$ | 1 | 1 | 1 | 1 | 0 |
| $\chi_1(n)$ | 1 | $-1$ | $-1$ | 1 | 0 |
| $\chi_2(n)$ | 1 | $i$ | $-i$ | $-1$ | 0 |
| $\chi_3(n)$ | 1 | $-i$ | $i$ | $-1$ | 0 |

Table 1.3: $q = 5$, $\phi(5) = 4$

| $n$ | 1 | 2 | 3 | 4 | 5 | 6 |
|---|---|---|---|---|---|---|
| $\chi_0(n)$ | 1 | 0 | 0 | 0 | 1 | 0 |
| $\chi_1(n)$ | 1 | 0 | 0 | 0 | $-1$ | 0 |

Table 1.4: $q = 6$, $\phi(6) = 2$

Not only can we describe Dirichlet characters as either principal or non-principal, there are other classifications depending on the character's specific properties. For example we can describe a Dirichlet character as either even or odd.

**Definition 1.7.** Let $\chi$ be a Dirichlet character modulo $q$.
We call $\chi$ **odd** if $\chi(-1) = -1$ or **even** if $\chi(-1) = 1$.

We can also define primitive Dirichlet characters, just as we defined primitive roots in §1.2. In the same way that a primitive root generates the coprime integers, all Dirichlet characters can be viewed as extensions of primitive Dirichlet characters.

**Definition 1.8.** Let $\chi$ be a Dirichlet character modulo $q$ and let $d \mid q$. Then $d$ is an **induced modulus** for $\chi$ if for all $a$ such that $(a, q) = 1$ and $a \equiv 1 \pmod{d}$ we have $\chi(a) = 1$.
A Dirichlet character is called **primitive** if it has no induced moduli. In other words, $\chi$ is primitive if and only if for all $d \mid q$ there exists an $a$ with $a \equiv 1 \pmod{d}$ and $(a, q) = 1$ such that $\chi(a) \neq 1$.

As we will see in later chapters, Dirichlet characters often appear to us in sums. We sometimes have to sum $\chi(n)$ over $n$ or perhaps sum over all the Dirichlet characters of the same order for a fixed $n$. We will see in the following part of this section some nice properties of the Dirichlet characters and their sums.

**Definition 1.9.** Let $\chi$ be any Dirichlet character modulo $m$ then

$$G(n, \chi) = \sum_{k=1}^{m} \chi(k) e^{2\pi i k n/m}$$

is called the **Gauss sum** associated with $\chi$.

The Gauss sum will be important in the next section as it is needed in the proof of the Pólya–Vinogradov inequality. The following proposition is important in proving the indicator function for primitive roots.

**Proposition 1.8.** *Let $\Gamma_d$ denote the set of all Dirichlet characters modulo $p$ of order $d$ and define*

$$S(d) = \sum_{\chi \in \Gamma_d} \chi(n).$$

*$S(d)$ is multiplicative.*

*Proof.* Let $d_1$ and $d_2$ be coprime integers and consider

$$S(d_1 d_2) = \sum_{\chi \in \Gamma_{d_1 d_2}} \chi(n) = \chi_0(n) + \chi_2(n) + \cdots + \chi_{\phi(d_1 d_2)-1}(n).$$

There are $\phi(d_1 d_2) = \phi(d_1)\phi(d_2)$ characters of order $d_1 d_2$ as $\phi$ is multiplicative (Proposition 1.3) and so there are $\phi(d_1)\phi(d_2)$ terms in the sum.

Now let $\psi_i$, where $i = 1, \ldots, \phi(d_1)$, be the $\phi(d_1)$ Dirichlet characters of order $d_1$ and let $\eta_j$, where $j = 1, \ldots, \phi(d_2)$, be the $\phi(d_2)$ Dirichlet characters of order $d_2$. Then

$$S(d_1)S(d_2) = \left( \sum_{\chi \in \Gamma_{d_1}} \chi(n) \right) \left( \sum_{\chi \in \Gamma_{d_2}} \chi(n) \right) = \sum_{i=1}^{\phi(d_1)} \sum_{j=1}^{\phi(d_2)} (\psi_i \eta_j)(n).$$

This sum has at most $\phi(d_1)\phi(d_2)$ terms. Therefore if we can show that the product of Dirichlet characters, $\psi\eta$, has order $d_1 d_2$ and that the sum $S(d_1)S(d_2)$ has exactly $\phi(d_1)\phi(d_2)$ terms, we have $S(d_1 d_2) = S(d_1)S(d_2)$.

Firstly we will show that if $\psi \in \Gamma_{d_1}$ and $\eta \in \Gamma_{d_2}$ then $\psi\eta \in \Gamma_{d_1 d_2}$.
Let $\psi \in \Gamma_{d_1}$ and $\eta \in \Gamma_{d_2}$, then $\psi^{d_1} = \chi_0$ and $\eta^{d_2} = \chi_0$ where $\chi_0$ is the principal character modulo $p$. So $(\psi\eta)^{d_1 d_2} = \psi^{d_1 d_2}\eta^{d_1 d_2} = \chi_0^{d_2}\chi_0^{d_1} = \chi_0$ which means $\mathrm{ord}(\psi\eta) \leq d_1 d_2$ (where $\mathrm{ord}(\chi)$ denotes the order of $\chi$). Suppose $\mathrm{ord}(\psi\eta) = K$ then by the division algorithm (Theorem 1.2 in [27]) $d_1 d_2 = qK + r$ where $q > 0$ and $0 \leq r < K$. Then

$$\chi_0 = (\psi\eta)^{d_1 d_2} = (\psi\eta)^{qK+r} = \chi_0(\psi\eta)^r = (\psi\eta)^r.$$

This means that if $1 \leq r < K$ then $\mathrm{ord}(\psi\eta) \leq r < K$. This is a contradiction as $\mathrm{ord}(\psi\eta) = K$. So $r = 0$, in particular $K \mid d_1 d_2$. So there exists $A$ such that $K = \dfrac{d_1 d_2}{A}$. Since the order is the least exponent, $K$, such that $(\psi\eta)^K = \chi_0$, $A$ is the greatest common divisor of $d_1$ and $d_2$. So $K = \dfrac{d_1 d_2}{(d_1,\, d_2)}$, namely $K$ is the least common multiple of $d_1$ and $d_2$. Since $(d_1, d_2) = 1$, $K = d_1 d_2$ and therefore $\psi\eta \in \Gamma_{d_1 d_2}$.

Now we will prove that the sum $S(d_1)S(d_2)$ has exactly $\phi(d_1)\phi(d_2)$ terms. Suppose $S(d_1)S(d_2)$ has less than $\phi(d_1)\phi(d_2)$ terms, in particular there is a double up of characters. Then $\psi_a\eta_b = \psi_c\eta_d$ defines a double up in the following three cases: ($a = c$ and $b \neq d$), ($a \neq c$ and $b = d$) and ($a \neq c$ and $b \neq d$). Without loss of generality we assume, $a \neq c$. Since $\psi_a\eta_b = \psi_c\eta_d$, $\psi_a = \psi_c\eta_d\eta_b^{-1}$ and we have two cases, where either $\eta_d\eta_b^{-1} = \chi_0$ or $\eta_d\eta_b^{-1} \neq \chi_0$.

Suppose $\eta_d\eta_b^{-1} = \chi_0$ then $\psi_a = \psi_c$ however this implies that $a = c$ which is a contradiction. Now suppose $\eta_d\eta_b^{-1} \neq \chi_0$ then since $\eta_b \in \Gamma_{d_2}$, $\eta_b^{-1} \in \Gamma_{d_2}$. So we have $\eta_d\eta_b^{-1} \in \Gamma_{d_2}$, that is $\text{ord}(\eta_d\eta_b^{-1}) = d_2$. If $d_2 = 1$ then $\eta_d\eta_b^{-1} - \chi_0$ and so $d_2 > 1$. Therefore by the first part of this proof $\text{ord}(\psi_c\eta_d\eta_b^{-1}) = \text{ord}(\psi_c)\text{ord}(\eta_d\eta_b^{-1}) = d_1 d_2 > d_1$. However $\psi_a = \psi_c\eta_d\eta_b^{-1}$ and $\text{ord}(\psi_a) = d_1$ which is a contradiction and so each $\psi_i\eta_j$ defines a unique character of order $d_1 d_2$ for all $1 \leq i \leq \phi(d_1)$ and $1 \leq j \leq \phi(d_2)$. Therefore the sum $S(d_1)S(d_2)$ has $\phi(d_1)\phi(d_2)$ unique terms.
Hence $S(d)$ is multiplicative. □

The following section will introduce a famous inequality on the sum of Dirichlet characters. This bound will play an important part in Chapter 2.

### 1.2.4 Pólya–Vinogradov inequality

The Pólya–Vinogradov inequality provides a bound on the sum of Dirichlet characters that is independent of the interval of summation. Because of this we can use this inequality when we need to bound the sum of Dirichlet characters in §2.1. There has been extensive research done on the inequality and therefore quite strong bounds are known. We will go through some of these in this section.

Let $\chi$ be a non-principal Dirichlet character modulo $q$ then the Pólya–Vinogradov inequality (Theorem 9.18 in [24]) states that

$$\sum_{n=M+1}^{M+N} \chi(n) = O\left(\sqrt{q}\log q\right) \tag{1.2}$$

for any integers $M$ and $N$ with $N > 0$. Equivalently,

$$\left|\sum_{n=M+1}^{M+N} \chi(n)\right| \leq c\sqrt{q}\log q \tag{1.3}$$

for a universal constant $c$.

This inequality was independently discovered by Pólya [28] and Vinogradov [32] in 1918. Consider the trivial bound on the same sum. Since $|\chi(n)| \leq 1$ for all $n$ (Proposition 1.6) we have

$$\left|\sum_{n=1}^{N} \chi(n)\right| = |\chi(1) + \chi(2) + \cdots + \chi(N)| \leq N. \tag{1.4}$$

We can expect the bound to be much smaller for non-principal characters. There is a lot of cancellation as the sum cycles through roots of unity. For example, recall the Dirichlet characters modulo 5 from Table 1.3. Let $N = 7$, then we have

$$\begin{aligned}
\sum_{n=1}^{7} \chi_0(n) &= \chi_0(1) + \chi_0(2) + \cdots + \chi_0(5) + \chi_0(1) + \chi_0(2) \\
&= 1 + 1 + 1 + 1 + 0 + 1 + 1 = 6, \\
\sum_{n=1}^{7} \chi_1(n) &= 1 - 1 - 1 + 1 + 0 + 1 - 1 = 0, \\
\sum_{n=1}^{7} \chi_2(n) &= 1 + i - i - 1 + 0 + 1 + i = 1 + i, \\
\sum_{n=1}^{7} \chi_3(n) &= 1 - i + i - 1 + 0 + 1 - i = 1 - i.
\end{aligned}$$

Note that if $N = 5$ then the sum equals 0 for all non-principal characters and $\phi(5) = 4$ for the principal character. Consider summing up to a multiple of the modulus. Let $\chi$ be a Dirichlet character modulo $q$, then for any $k \geq 1$, by periodicity we have

$$\sum_{n=1}^{kq} \chi(n) = \begin{cases} 0 & \text{if } \chi \text{ is non-principal,} \\ k\phi(q) & \text{if } \chi \text{ is principal.} \end{cases}$$

We can see from the above example that when $\chi$ is non-principal there is a lot cancellation however when $\chi$ is a principal character we are summing a string of ones, with zeros appearing only when $(n, q) > 1$. Therefore we only need to bound the sum when $\chi$ is non-principal as it is known when $\chi$ is principal.

As a result of Montgomery and Vaughan's work in [23] we have the following lower bound on the character sum

$$\left|\sum_{n=M+1}^{M+N} \chi(n)\right| \gg \sqrt{q}$$

for any primitive character $\chi$ modulo $q$. This shows that apart from the factor of log, the inequality (1.2) is the best possible. Assuming the Generalised Riemann

Hypothesis[1], Montgomery and Vaughan [22] have also shown that

$$\left|\sum_{n=M+1}^{M+N} \chi(n)\right| \ll \sqrt{q}\log\log q$$

for any Dirichlet character modulo $q$. The implicit constant for the Pólya–Vinogradov inequality can be shown to be 1 for all non-principal characters. One proof of this can be found in Davenport's book [7]. We will prove $c = 1$ for all primitive characters as the full proof for all non-principal characters is long and can be extended from the following proof.

**Theorem 1.1** (Theorem 8.12 in [1]). *Let $\chi$ be any primitive Dirichlet character modulo $q$ then for all $x \geq 1$ we have*

$$\left|\sum_{n\leq x} \chi(n)\right| < \sqrt{q}\log q.$$

*Proof.* Let $\chi$ be a primitive character modulo $q$ then the finite Fourier expansion of $\chi(n)$ (Theorem 8.20 in [1]) is

$$\chi(n) = \frac{G(1,\,\chi)}{q}\sum_{k=1}^{q}\bar{\chi}(k)e^{-2\pi ink/q}$$

where the Gauss sum is

$$G(1,\,\chi) = \sum_{n=1}^{q}\chi(n)e^{2\pi in/q}$$

and $\bar{\chi}$ is the complex conjugate of $\chi$.

Now we sum over all $n \leq x$ to obtain

$$\sum_{n\leq x}\chi(n) = \frac{G(1,\,\chi)}{q}\sum_{k=1}^{q-1}\bar{\chi}(k)\sum_{n\leq x}e^{-2\pi ink/q}.$$

Here the sum of Dirichlet character $\chi(k)$ is now between $1 \leq k \leq q-1$ since $\chi(q) = 0$. Since we are looking for a bound we need to take absolute values, which results in

$$\begin{aligned}\left|\sum_{n\leq x}\chi(n)\right| &= \left|\frac{G(1,\,\chi)}{q}\sum_{k=1}^{q-1}\bar{\chi}(k)\sum_{n\leq x}e^{-2\pi ink/q}\right| \\ &\leq \frac{|G(1,\,\chi)|}{q}\sum_{k=1}^{q-1}\left|\bar{\chi}(k)\sum_{n\leq x}e^{-2\pi ink/q}\right|.\end{aligned}$$

[1] In 1859 Riemann conjectured that the nontrivial zeros of the zeta function have real part 1/2. The Generalised Riemann Hypothesis concerns the zeros of L-functions, which are similar to the zeta function[24].

From Proposition 1.6 we have $|\chi(n)| \leq 1$ and so it follows that $|\bar{\chi}(n)| \leq 1$ and so

$$\left|\sum_{n\leq x}\chi(n)\right| \leq \frac{|G(1,\,\chi)|}{q}\sum_{k=1}^{q-1}\left|\sum_{n\leq x}e^{-2\pi ink/q}\right|$$
$$= \frac{|G(1,\,\chi)|}{q}\sum_{k=1}^{q-1}|f(k)| \tag{1.5}$$

where

$$f(k) = \sum_{n\leq x}e^{-2\pi ink/q}.$$

Now consider the function $f(k)$.

$$f(q-k) = \sum_{n\leq x}e^{-2\pi in(q-k)/q} = \sum_{n\leq x}e^{2\pi ink/q}e^{-2\pi in} = \sum_{n\leq x}e^{2\pi ink/q} = \overline{f(n)}$$

which means $|f(q-k)| = |\overline{f(k)}| = |f(k)|$. Hence

$$\sum_{k=1}^{q-1}|f(k)| \leq \sum_{k=1}^{q/2}|f(k)| + \sum_{k=q/2}^{q-1}|f(q-k)|$$
$$= 2\sum_{k=1}^{q/2}|f(k)|. \tag{1.6}$$

Now let $r = e^{-2\pi ik/q}$ and $a = [x]$ then $f(k)$ is a geometric series,

$$f(k) = \sum_{n=1}^{a}r^n$$

with $r \neq 1$ since $1 \leq k \leq q-1$ and writing $t = e^{-\pi ik/q}$ we have $r = t^2$ and $t^2 \neq 1$ since $1 \leq k \leq q/2$. Hence

$$f(k) = \frac{r(1-r^a)}{1-r} = \frac{t^2(1-t^{2a})}{1-t^2} = \frac{t^{2+a}(t^{-a}-t^a)}{t(t^{-1}-t)} = t^{1+a}\frac{t^{-a}-t^a}{t^{-1}-t}$$

and so using Euler's formula, $e^{ix} = \cos x + i\sin x$, we have

$$|f(k)| = \left|\frac{t^{-a}-t^a}{t^{-1}-t}\right| = \left|\frac{e^{\pi ika/q}-e^{-\pi ika/q}}{e^{\pi ik/q}-e^{-\pi ik/q}}\right| = \frac{\left|\sin\left(\frac{\pi ka}{q}\right)\right|}{\left|\sin\left(\frac{\pi k}{q}\right)\right|} \leq \frac{1}{\sin\left(\frac{\pi k}{q}\right)}. \tag{1.7}$$

Now we relax the bound on $|f(k)|$ by using the inequality $\sin(u) \geq 2u/\pi$ which is valid for $0 \leq u \leq \pi/2$. Since $k \leq q/2$ then $\pi k/q \leq \pi/2$ and so (1.7) becomes

$$|f(k)| \leq \frac{1}{\frac{2}{\pi}\frac{\pi k}{q}} = \frac{q}{2k}. \tag{1.8}$$

Substituting (1.6) and (1.8) into (1.5) we have

$$\begin{aligned}\left|\sum_{n\le x}\chi(n)\right| &\le \frac{|G(1,\,\chi)|}{q}2\sum_{k=1}^{q/2}\frac{q}{2k}\\ &= |G(1,\,\chi)|\sum_{k=1}^{q/2}\frac{1}{k}\\ &< |G(1,\,\chi)|\log q. \end{aligned} \tag{1.9}$$

Since $\chi$ is a primitive Dirichlet character modulo $q$, for all $n$ such that $(n,\,q)=1$ we have $|G(n,\,\chi)|=\sqrt{q}$. A proof of this can be found in Theorem 8.11 and Theorem 8.19 from [1]. Therefore

$$\left|\sum_{n\le x}\chi(n)\right| < \sqrt{q}\log q.$$

□

There has been extensive work done on improving the upper bound of (1.3). These improvements are obtained from advanced methods and so they will not be proved here. In 2007 Granville and Soundararajan [11] showed that if $\chi$ has odd order then a small power of $\log q$ can be saved in (1.3). The following result was improved in 2012 by Goldmakher [10] who stated that for each fixed odd number $g>1$, for $\chi$ of order $g$,

$$\left|\sum_{n\le x}\chi(n)\right| \le \sqrt{q}(\log q)^{\Delta_g+O(1)} \qquad \text{where } \Delta_g=\frac{g}{\pi}\sin\frac{\pi}{g}, \text{ as } q\to\infty.$$

In 2013, Frolenkov and Soundararajan [9] were able to obtain an explicit version of the Pólya–Vinogradov inequality for all non-principal Dirichlet characters. They used a result which bounds the sum (over any interval $[M+1,\,M+N]$) of a broader class of arithmetic functions than the Dirichlet characters to obtain the following theorem.

**Theorem 1.2** (Corollary 1 in [9]). *Let $q>100$ and let $\chi$ be any non-principal Dirichlet character modulo $q$. Then we have*

$$\left|\sum_{n\le x}\chi(n)\right| < c\sqrt{q}\log q$$

*where*

$$c=\left(\frac{1}{\pi\sqrt{2}}+\frac{6}{\pi\sqrt{2}\log q}+\frac{1}{\log q}\right). \tag{1.10}$$

The parameter, $c$ will decrease as we take $q$ to be large. In particular, both the second and the third term tend to zero and $c \to (\pi\sqrt{2})^{-1}$ as $q \to \infty$. This will be important when we use this theorem later in §2.1.

Recall the trivial bound (1.4). There are some cases where the trivial bound is lower than the bound in Theorem 1.2. Hence, we write

$$\left| \sum_{n \le x} \chi(n) \right| < \min \{x,\, c\sqrt{q}\log q\}.$$

There have been many improvements to (1.3) for primitive Dirichlet characters. In particular, sharper bounds have been found when splitting the sum into the sum of either even or odd Dirichlet characters. Explicit bounds of this form have been found by Pomerance [29] and Frolenkov [8] in 2011. As a result of the same method used to prove Theorem 1.2, Frolenkov and Soundararajan proved the following theorem.

**Theorem 1.3** (Theorem 2 in [9]). *Let $\chi$ be a primitive Dirichlet character modulo $q$. If $\chi$ is even and $q \ge 1200$ we have*

$$\left| \sum_{n \le x} \chi(n) \right| \le \frac{2}{\pi^2}\sqrt{q}\log q + \sqrt{q}.$$

*If $\chi$ is odd and $q \ge 40$ we have*

$$\left| \sum_{n \le x} \chi(n) \right| \le \frac{1}{2\pi}\sqrt{q}\log q + \sqrt{q}.$$

It should be noted that if $q$ is taken to be much larger, say greater $10^6$, then there are some mild improvements to the constants in Theorem 1.2 and Theorem 1.3. However these improvements have almost no effect on the results of this thesis.

Now that we have introduced all the necessary background information we are able to discuss the distribution of primitive roots of primes in the next chapter.

# 2

# Least square-free primitive root

There are many questions concerning primitive roots. There is a famous conjecture by Artin [13] that states that given an integer that is neither $-1$ nor a perfect square, that integer is a primitive root of infinitely many primes. Heath-Brown [16] proved in 1985 that Artin's conjecture fails for at most two primes $p$. Heath-Brown also proved that there are at most three square-free integers for which Artin's conjecture fails. Pieter Moree provides a survey [25] of the results on Artin's conjecture.

The distribution of primitive roots modulo a prime is also of interest and this is where the results of this thesis fit in. In particular we will be looking at the least primitive root modulo a prime. There is a well known conjecture from Erdős [13]. He asks: do all primes $p$ have a prime primitive root less than $p$? This conjecture is one of the many unsolved problems of primitive roots. Let $g(p)$ denote the least primitive root modulo prime $p$. Numerical examples [26] show that we expect $g(p)$ to be very small. For example among the first $19,862$ primes, 37.4% of these primes have $g(p) = 2$ and 22.8% of the primes have $g(p) = 3$. We actually have that 80% of the first $19,862$ primes have $g(p) \leq 6$. In 1961 Burgess [2] proved that for any fixed

$\epsilon > 0$ we have

$$g(p) = O\left(p^{1/4+\epsilon}\right).$$

That is to say that for sufficiently large $p$, the least primitive root of $p$ is very small. More recently an asymptotic bound for the least prime primitive root modulo $p$ has also been found. Let $\hat{g}(p)$ denote the least prime primitive root modulo prime $p$, then in 2015 Ha [14] proved

$$\hat{g}(p) = O\left(p^{3.1}\right).$$

This bound does not tell us much about $\hat{g}(p)$ because for large primes this bound is huge. We also see from this bound that we are a long way off from solving Erdős' conjecture. We expect the asymptotic bound on $\hat{g}(p)$ to be small, assuming the Generalised Riemann Hypothesis, Shoup proved that

$$\hat{g}(p) = O\left((\log p)^6\right).$$

Although we expect the least prime primitive root modulo a prime to be small, it is very difficult to improve the asymptotic bound. However there have been some explicit improvements. Grosswald conjectured in 1981 [12] that

$$g(p) < \sqrt{p} - 2 \quad \text{for all } p > 409.$$

There has been some recent work on resolving Grosswald's conjecture. In 2016 Cohen, Oliveira e Silva and Trudgian proved Grosswald's conjecture for $409 < p < 2.5 \times 10^{15}$ and $p > 3.38 \times 10^{71}$ [6]. They also prove the following bound for the least prime primitive root modulo a prime.

**Theorem 2.1** (§4 of [6]). *Given a prime $p$ we have*

$$\hat{g}(p) < \sqrt{p} - 2 \quad \textit{for} \quad 2791 < p < 2.5 \times 10^{15}.$$

Assuming the Generalised Riemann Hypothesis, McGown, Treviño and Trudgian [21] proved $\hat{g}(p) < \sqrt{p} - 2$ for all primes $p > 2791$ and $g(p) < \sqrt{p} - 2$ for all primes $p > 409$. That is, Grosswald's conjecture is true assuming the Generalised Riemann Hypothesis.

It is a difficult problem to obtain a bound on the least prime primitive root modulo a prime and that is why we consider square-free primitive roots of primes. These can be seen as a slight generalisation of prime primitive roots and there is room for improvement on their bound.

Let $g^{\square}(p)$ denote the least square-free primitive root modulo prime $p$. Just like the least prime primitive root there has been some work on bounding $g^{\square}(p)$ explicitly and implicitly. Shapiro [31] proved that

$$g^{\square}(p) = O\left(p^{1/2+\epsilon}\right).$$

This bound was improved slightly in 2005 by Liu and Zhang [19]. They showed

$$g^{\square}(p) = O\left(p^{9/22+\epsilon}\right)$$

and Martin [20] has suggested that the bound on the least square-free primitive root is of the same form as Burgess bound on $g(p)$. That is $g^{\square}(p) = O\left(p^{1/4+\epsilon}\right)$ for some fixed $\epsilon > 0$. Recently there have been results published on an explicit bound for the least square-free primitive root modulo a prime. Cohen and Trudgian (Theorem 1 in [6]) proved that

$$g^{\square}(p) < p^{0.96} \quad \text{for all primes } p.$$

This means that all primes have a square-free primitive root less than themselves. One question that arises from this bound is, can we achieve a lower exponent? This question is what motivates the next two chapters. In the following two chapters we will outline how the following question can be answered.

**Question 2.1.** For what $\alpha < 1$ is the following statement true?

All primes $p$ have a square-free primitive root less than $p^{\alpha}$.

It should be noted that there is a lower bound for $\alpha$. Since 2 is the least square-free primitive root of 3, the above theorem states $2 < 3^{\alpha}$ and so we have $\alpha > \log 2/\log 3 > 0.6309$.

To answer Question 2.1 for a given $\alpha$ we follow a similar structure to other proofs on bounds on least primitive roots of primes, for example Cohen and Trudgian's paper [6] and McGown, Treviño and Trudgian's paper [21]. These proofs have three steps. We will outline the first step in §2.1, the second and third steps will be outlined in §3.2 and §3.3 respectively. At each stage there is an improvement made to the $\alpha = 0.96$ obtained in [6].

## 2.1 Results without sieving

In this section we will outline the first step in the proof of Question 2.1. We first find an expression for $N^{\square}(p, x)$ defined to be the number of square-free primitive roots modulo $p$ less than $x$. To achieve this we need to find an indicator function for primitive roots of primes.

**Lemma 2.1** (Lemma 8.5.1 from [31]). *Let $p$ be an odd prime and let $\Gamma_d$ denote the set of Dirichlet characters modulo $p$ of order $d$. Then*

$$\frac{\phi(p-1)}{p-1}\sum_{d|p-1}\frac{\mu(d)}{\phi(d)}\sum_{\chi\in\Gamma_d}\chi(n)=\begin{cases}1 & \text{if } n \text{ is a primitive root modulo } p,\\ 0 & \text{otherwise.}\end{cases}$$

*Proof.* Let $q_1$, $q_2, \ldots, q_r$ be the distinct prime divisors of $p-1$. Note both $\mu(d)$ and $\phi(d)$ are multiplicative by Proposition 1.3. Also by Proposition 1.8 $\sum_{\chi\in\Gamma_d}\chi(n)$ is multiplicative in $d$.
Therefore $\sum_{d|p-1}\frac{\mu(d)}{\phi(d)}\sum_{\chi\in\Gamma_d}\chi(n)$ is multiplicative in $p-1$ and we have

$$\sum_{d|p-1}\frac{\mu(d)}{\phi(d)}\sum_{\chi\in\Gamma_d}\chi(n)=\left(\sum_{d|q_1^{\alpha_1}}\frac{\mu(d)}{\phi(d)}\sum_{\chi\in\Gamma_d}\chi(n)\right)\cdots\left(\sum_{d|q_r^{\alpha_r}}\frac{\mu(d)}{\phi(d)}\sum_{\chi\in\Gamma_d}\chi(n)\right)$$

where $\alpha_i\geq 1$. Since $\mu(q^\alpha)=0$ for all $\alpha\geq 2$, we are left with

$$\sum_{d|p-1}\frac{\mu(d)}{\phi(d)}\sum_{\chi\in\Gamma_d}\chi(n)=\prod_{j=1}^{r}\left(\frac{\mu(1)}{\phi(1)}\sum_{\chi\in\Gamma_1}\chi(n)+\frac{\mu(q_j)}{\phi(q_j)}\sum_{\chi\in\Gamma_{q_j}}\chi(n)\right)$$

Recall the only character of order 1 is the principal character and as the modulus is prime we have $\chi_0(n)=1$. Also recall $\mu(p)=-1$ for all primes $p$ and $\mu(1)=1=\phi(1)$. Therefore

$$\sum_{d|p-1}\frac{\mu(d)}{\phi(d)}\sum_{\chi\in\Gamma_d}\chi(n)=\prod_{j=1}^{r}\left(1-\frac{1}{\phi(q_j)}\sum_{\chi\in\Gamma_{q_j}}\chi(n)\right). \tag{2.1}$$

Now fix $g$ a primitive root of $p$ and let $a=\text{ind}_g n$ and $\chi(n)\in\Gamma_d$. Then as $\chi$ is multiplicative $\chi(n)=\chi(g^a)=(\chi(g))^a=e^{2\pi ika/d}$ for $1\leq k\leq d$ such that $(k,\, d)=1$, and

$$\frac{1}{\phi(q_j)}\sum_{\chi\in\Gamma_{q_j}}\chi(n)=\frac{1}{q_j-1}\sum_{k=1}^{q_j-1}e^{2\pi ika/q_j}.$$

Suppose $q_j\mid a$, then $e^{2\pi ika/q_j}=1$ and so $\sum_{k=1}^{p-1}e^{2\pi ika/q_j}=q_j-1$. Now suppose $q_j\nmid a$, then since $k<q_j-1$, $q_j\nmid ka$ and $e^{2\pi ika/q_j}\neq 1$ for all $k$.

Consider $\sum_{k=1}^{q_j-1} X^k$ where $X = e^{2\pi i k a/q_j}$. This is a geometric series and so we have

$$\sum_{k=1}^{q_j-1} X^k = \frac{X - X^{q-j}}{1-X} = \frac{X-1}{1-X} = -1.$$

Therefore

$$\frac{1}{q_j-1}\sum_{k=1}^{p-1} e^{2\pi i k a/q_j} = \begin{cases} 1 & \text{if } q_j \mid a, \\ -\frac{1}{q_j-1} & \text{if } q_j \nmid a, \end{cases}$$

which gives us

$$1 - \frac{1}{\phi(q_j)} \sum_{\chi \in \Gamma_{q_j}} \chi(n) = \begin{cases} 0 & \text{if } q_j \mid \text{ind}_g n, \\ \frac{q_j}{q_j-1} & \text{if } q_j \nmid \text{ind}_g n. \end{cases} \tag{2.2}$$

Suppose $q_j \nmid \text{ind}_g n$ for all $j$, then by Proposition 1.5

$$\prod_{j=1}^{r} \left(1 - \frac{1}{\phi(q_j)} \sum_{\chi \in \Gamma_{q_j}} \chi(n)\right) = \prod_{j=1}^{r} \frac{q_j}{q_j - 1} = \frac{p-1}{\phi(p-1)}$$

which combined with 2.2 and 2.1 gives us

$$\frac{\phi(p-1)}{p-1} \sum_{d \mid p-1} \frac{\mu(d)}{\phi(d)} \sum_{\chi \in \Gamma_d} \chi(n) = \begin{cases} 1 & \text{when } q_j \nmid \text{ind}_g n \quad \text{for all} \quad j, \\ 0 & \text{when } q_j \mid \text{ind}_g n \text{ for some } j. \end{cases}$$

The condition $q_j \nmid \text{ind}_g n$ for all $j$ implies $(\text{ind}_g n,\, p-1) = 1$ which is the condition for $n$ to be a primitive root modulo $p$ (Remark 1.1). Hence

$$\frac{\phi(p-1)}{p-1} \sum_{d \mid p-1} \frac{\mu(d)}{\phi(d)} \sum_{\chi \in \Gamma_d} \chi(n) = \begin{cases} 1 & \text{if } n \text{ is a primitive root (mod } p), \\ 0 & \text{otherwise.} \end{cases}$$

□

Now using this indicator function for primitive roots, we can sum over the square-free integers to obtain

$$\begin{aligned} N^{\square}(p,\, x) &= \sum_{\substack{n \le x \\ n = \square-\text{free}}} f(n) \\ &= \frac{\phi(p-1)}{p-1} \left\{ \sum_{\substack{n \le x \\ n = \square-\text{free}}} 1 + \sum_{\substack{d \mid p-1 \\ d > 1}} \frac{\mu(d)}{\phi(d)} \sum_{\chi \in \Gamma_d} \sum_{\substack{n \le x \\ n = \square-\text{free}}} \chi(n) \right\}. \end{aligned} \tag{2.3}$$

Therefore to prove Question 2.1 for a given $\alpha$ we need $N^{\square}(p, x) > 0$ for $x = p^{\alpha}$ and from (2.3) we have

$$N^{\square}(p, x) > 0 \iff \sum_{\substack{n \le x \\ n = \square - \text{free}}} 1 + \sum_{\substack{d | p-1 \\ d > 1}} \frac{\mu(d)}{\phi(d)} \sum_{\chi \in \Gamma_d} \sum_{\substack{n \le x \\ n = \square - \text{free}}} \chi(n) > 0. \tag{2.4}$$

Now we need to find a bound on the right hand side of (2.4). First we use (1.1) to separate the innermost sum and then reverse summation to obtain

$$\sum_{\substack{n \le x \\ n = \square - \text{free}}} \chi(n) = \sum_{n \le x} \chi(n) \sum_{d^2 | n} \mu(d) = \sum_{1 \le d \le \sqrt{x}} \mu(d) \sum_{\substack{n \le x \\ n \equiv 0 \pmod{d^2}}} \chi(n). \tag{2.5}$$

Consider the trivial bound on the innermost sum. From Proposition 1.6 we have that $|\chi(n)| \le 1$ for all $n$. This bound together with the multiplicativity of $\chi$ means

$$\begin{aligned} \left| \sum_{\substack{n \le x \\ n \equiv 0 \pmod{d^2}}} \chi(n) \right| &= \left| \chi(d^2) + \chi(2d^2) + \cdots + \chi\left( \left[ \frac{x}{d^2} \right] d^2 \right) \right| \\ &\le |\chi(d^2)| + |\chi(2d^2)| + \cdots + \left| \chi\left( \left[ \frac{x}{d^2} \right] d^2 \right) \right| \\ &\le \frac{x}{d^2}. \end{aligned}$$

Recall Theorem 1.2, a version of the Pólya–Vinogradov inequality which also provides an upper bound on this sum. Then for $p > 100$ we have

$$\left| \sum_{\substack{n \le x \\ n \equiv 0 \pmod{d^2}}} \chi(n) \right| \le \min \left\{ \frac{x}{d^2}, \, c\sqrt{p} \log p \right\}$$

where

$$c = \left( \frac{1}{\pi\sqrt{2}} + \frac{6}{\pi\sqrt{2}\log p} + \frac{1}{\log p} \right).$$

Note that it is fine that we are taking $p > 100$ here as we will show later that we will be using this bound for $p > 2.5 \times 10^{15}$.

Now we can use this bound to separate (2.5) into two parts, choosing to sum $x/d^2$ over $d > d_0$ and $c\sqrt{p}\log p$ over $d \le d_0$. Here $d_0$ is chosen to obtain the smallest bound. Separating the sum we obtain

$$\left| \sum_{\substack{n \le x \\ n = \square - \text{free}}} \chi(n) \right| \le \sum_{d \le d_0} |\mu(d)| c\sqrt{p} \log p + \sum_{d_0 < d \le \sqrt{x}} |\mu(d)| \frac{x}{d^2}. \tag{2.6}$$

We can estimate the first sum of (2.6) using Cipu's result ((*a*) from Lemma 1.1),

$$\sum_{d\le d_0}|\mu(d)|c\sqrt{p}\log p \le \left(\frac{6}{\pi^2}+A\sqrt{d_0}\right)c\sqrt{p}\log p \tag{2.7}$$

where $A = 0.679091$. We keep the equation in terms of the general constant $A$ so that we can investigate whether the stronger bounds for larger $x$ (from Lemma 1.1) will make a significant difference to our results. We will discuss this at the end of the chapter.

Now using partial summation and Cipu's result we can estimate the second sum.

$$\begin{aligned}\sum_{d_0<d\le\sqrt{x}}|\mu(d)|\frac{x}{d^2} &= x\left(\frac{1}{x}\sum_{d\le\sqrt{x}}|\mu(d)|-\frac{1}{d_0^2}\sum_{d\le d_0}|\mu(d)|+2\int_{d_0}^{\sqrt{x}}\left(\sum_{n\le t}|\mu(n)|\right)t^{-3}dt\right)\\ &\le\left(\frac{6}{\pi^2}\sqrt{x}+Ax^{1/4}\right)-\left(\frac{6}{\pi^2}d_0-A\sqrt{d_0}\right)\frac{x}{d_0^2}\\ &\quad+2x\int_{d_0}^{\sqrt{x}}\left(\frac{6}{\pi^2}t+A\sqrt{t}\right)t^{-3}dt.\end{aligned} \tag{2.8}$$

Integrating by parts, the last term of the right hand side of (2.8) becomes

$$\begin{aligned}2x\int_{d_0}^{\sqrt{x}}\left(\frac{6}{\pi^2}t+A\sqrt{t}\right)t^{-3}dt &= \frac{12}{\pi^2}x\int_{d_0}^{\sqrt{x}}t^{-2}dt+2Ax\int_{d_0}^{\sqrt{x}}t^{-5/2}dt\\ &= -\frac{12}{\pi^2}x\left(\frac{1}{\sqrt{x}}-\frac{1}{d_0}\right)-\frac{4}{3}Ax\left(x^{-3/4}-d_0^{-3/2}\right).\end{aligned} \tag{2.9}$$

Now substituting (2.9) into (2.8) the bound on the second sum of (2.6) is

$$\sum_{d_0<d\le\sqrt{x}}|\mu(d)|\frac{x}{d^2}\le-\frac{6}{\pi^2}\sqrt{x}+\frac{6}{\pi^2}\frac{x}{d_0}-\frac{1}{3}Ax^{1/4}+\frac{7}{3}Axd_0^{-3/2}. \tag{2.10}$$

Adding (2.7) to (2.10) we obtain the following bound

$$\left|\sum_{\substack{n\le x\\ n=\square-\text{free}}}\chi(n)\right|\le\left(\frac{6}{\pi^2}+A\sqrt{d_0}\right)c\sqrt{p}\log p+\frac{6}{\pi^2}\left(\frac{x}{d_0}-\sqrt{x}\right)-\frac{1}{3}Ax^{1/4}+\frac{7}{3}Axd_0^{-3/2}. \tag{2.11}$$

Now that we have a bound we need to find a point to separate the interval of summation, $d_0$, that is close to optimal and an integer. We would like to choose $d_0$ such that the bound (2.11) is minimised. This is approximately when $x/d_0^2 = c\sqrt{p}\log p$. Since our $d_0$ needs to be an integer we will take the integer part,

$$d_0=\left[\left(\frac{x}{c\sqrt{p}\log p}\right)^{1/2}\right].$$

Now let $D = \left(\frac{x}{c\sqrt{p}\log p}\right)^{1/2}$ then $D - 1 < [D] \leq D$ and (2.11) becomes

$$\begin{aligned}\left|\sum_{\substack{n\leq x\\ n=\square-\text{free}}} \chi(n)\right| &\leq \frac{6}{\pi^2}[D]c\sqrt{p}\log p + \frac{6}{\pi^2}\frac{x}{[D]} - \frac{6}{\pi^2}\sqrt{x} \\ &\quad + A\sqrt{[D]}c\sqrt{p}\log p - \frac{1}{3}Ac^{1/4} + \frac{7}{3}Ax[D]^{-3/2} \\ &\leq \frac{6}{\pi^2}Dc\sqrt{p}\log p + \frac{6}{\pi^2}\frac{x}{D-1} - \frac{6}{\pi^2}\sqrt{x} + A\sqrt{D}c\sqrt{p}\log p \\ &\quad - \frac{1}{3}Ax^{1/4} + \frac{7}{3}Ax(D-1)^{-3/2}. \end{aligned} \tag{2.12}$$

Let $\hat{D} = \frac{D}{D-1}$ so that the main term does not contain $D$ then (2.12) becomes

$$\begin{aligned}\left|\sum_{\substack{n\leq x\\ n=\square-\text{free}}} \chi(n)\right| &\leq \frac{6}{\pi^2}\sqrt{x}(c\sqrt{p}\log p)^{1/2} + \frac{6}{\pi^2}\hat{D}\sqrt{x}(c\sqrt{p}\log p)^{1/2} - \frac{6}{\pi^2}\sqrt{x} \\ &\quad + Ax^{1/4}(c\sqrt{p}\log p)^{3/4} - \frac{1}{3}Ax^{1/4} + \frac{7}{3}Ax\left(\left(\frac{x}{c\sqrt{p}\log p}\right)^{1/2} - 1\right)^{-3/2}. \end{aligned} \tag{2.13}$$

Recall from Theorem 1.2 that $c$ tends to $(\sqrt{2}\pi)^{-1}$ as $p$ increases. We also have that for large $p$, $D$ is large and so $\hat{D}$ tends to 1. When Cohen and Trudgian estimated this bound, in [6], they used the trivial bound for the Möbius function, $|\mu(d)| \leq 1$ in (2.6) and the constant $c = 1$ for the Pólya–Vinogradov inequality. Therefore we have made an improvement to the bound (2.13) which should help in lowering the exponent $\alpha$ in Question 2.1.

The above estimation for the innermost sum of (2.4) does not depend on Dirichlet characters. Therefore when substituting (2.13) into (2.4) we are summing a constant over all Dirichlet characters modulo $p$ of order $d$. Recall from Proposition 1.7 there are $\phi(d)$ Dirichlet characters modulo $p$ of order $d$. Therefore we get cancellation of $\phi(d)^{-1}$.

The first sum of the right hand side of (2.4) can be estimated using the lower bound from Cipu ($(a)$ from Lemma 1.1). We obtain

$$\sum_{\substack{n\leq x\\ n=\square-\text{free}}} 1 \geq \frac{6}{\pi^2}x - 0.104\sqrt{x}. \tag{2.14}$$

Again this bound can be increased by taking $x$ to be larger however we will discuss at the end of the chapter how the slight increase in the constant does not make a significant difference to our results. Hence we use the bound stated above as it holds for all $x$. The only sum in the right hand side of (2.4) left to estimate is

$$\sum_{\substack{d|p-1\\d>1}} |\mu(d)| = \text{ the number of square-free divisors of } p-1 \text{ excluding 1.}$$

Let $\omega(n)$ denote the number of distinct prime divisors of $n$, then $p-1$ has a prime decomposition with $\omega(p-1)$ distinct primes. Therefore a square-free divisor of $p-1$ will have a prime decomposition of the same distinct primes, each with an exponent of either 0 or 1. This results in there being $2^{\omega(p-1)}$ square-free divisors of $p-1$. Hence

$$\sum_{\substack{d|p-1\\d>1}} |\mu(d)| = 2^{\omega(p-1)} - 1. \tag{2.15}$$

Now substituting (2.14), (2.15) and (2.13) into (2.4) and setting $x = p^\alpha$ we have $N^\square(p, x) > 0$ if

$$G(x) := x^{1/2}p^{-1/4} - \frac{\pi^2}{6}\left(\frac{0.104}{p^{1/4}} + \left(2^{\omega(p-1)} - 1\right)E\right) > 0 \tag{2.16}$$

where

$$\begin{aligned} E = (\log p)^{1/2}\Bigg(&\frac{6}{\pi^2}(c^{1/2} + \hat{D}c^{1/2} - p^{-1/4}(\log p)^{-1/2}) \\ &+ A\Big(p^{1/8-1/4\alpha}c^{3/4}(\log p)^{1/4} - \frac{1}{3}p^{-1/4-1/4\alpha}(\log p)^{1/2} \\ &+ \frac{7}{3}p^{1/8-1/4\alpha}(\log p)^{1/4}c^{3/4}\Big)\Bigg) \qquad (2.17) \\ \le (\log p)^{1/2}\Bigg(&\frac{6}{\pi^2}c^{1/2}\ (1+\hat{D}) + \frac{10}{3}Ap^{1/8-\alpha/4}c^{3/4}(\log p)^{1/4}\Bigg). \end{aligned}$$

We can see that $G(x)$ has a main term and an error term and so to prove Question 2.1 for a given $\alpha$ we need to show that the main term outweighs the error term for all primes $p$, setting $x = p^\alpha$. Since $\alpha > 1/2$ we have that $E$ is equal to $(\log p)^{1/2}$ multiplied by a decreasing function in $p$. The function $G(x)$ is of the same form as (8) from Cohen and Trudgian's paper [6]. They obtained

$$G_{CT}(x) := x^{1/2}p^{-1/4} - \frac{\pi^2}{6}\left(\frac{0.104}{p^{1/4}} + 2^{\omega(p-1)+1}(\log p)^{1/2}\right).$$

We can see that in $E$ we have $6/\pi^2 c^{1/2}$ whereas in $G_{CT}$ they have 2. This results in a smaller error term in $G(x)$.

Note that for all computations, the stricter equation for $E$ (2.17) is used instead of the upper bound.

Now by substituting $x = p^\alpha$ into $G(x)$ we will be able to find a value $n = \omega(p-1)$ such that Question 2.1 is true for all $p$ such that $\omega(p-1) \geq n$ for that given $\alpha$. For example take $\alpha = 0.96$, then we have that $G(p^{0.96}) > 0$ for all $p$ such that $\omega(p-1) \geq 26$ (this is an improvement on 30 from [6]). Recall Theorem 2.1 which states that

$$\hat{g}(p) < \sqrt{p} - 2 \text{ for all primes } p \text{ satisfying } 2791 \leq p \leq 2.5 \times 10^{15}.$$

Now this implies that Question 2.1 with $\alpha > 0.6309$ is true for this interval of primes. Running a quick computation we find that all primes $p$ less than or equal to 2791 have a square-free primitive root less than $p^{0.6309}$. Hence we have that Question 2.1 is true for all primes less than $2.5 \times 10^{15}$. For $\alpha = 0.96$ this takes care of cases $1 \leq \omega(p-1) \leq 9$ (an improvement from [6]). For example when $\omega(p-1) = 9$, we have that $G(p^{0.96}) > 0$ for all $p > 2.48 \times 10^{15}$. This means that we only need to consider $p < 2.48 \times 10^{15}$ which are all covered by Theorem 2.1. We dispatch all cases $1 \leq \omega(p-1) \leq 9$ in this way. We have that all possible exceptions to Question 2.1 for $\alpha = 0.96$ occur when $10 \leq \omega(p-1) \leq 25$.

Table 2.1 shows for what values of $\omega(p-1)$ Question 2.1 was unable to be answered at this stage in the proof. As mentioned earlier in the section, there are better bounds from Cipu (Lemma 1.1) that can be used in the formulation of $G(x)$ however these stronger bounds do not change the intervals in Table 2.1 and therefore do not make a significant difference to our results. This is not surprising as, the explicit constants from Cipu's bound that we have used are already so small that they are insignificant in $G(x)$ compared to $p$ which is large. For example we have the explicit constant 0.104 which is divided by $p^{1/4}$ in $G(x)$. Now taking $p$ to be the smallest possible, that is $p = 2.5 \times 10^{15}$ we have

$$\frac{0.104}{p^{1/4}} = 0.0000147$$

and if we replace 0.104 with the smallest explicit bound from Cipu we have

$$\frac{0.02767}{p^{1/4}} = 0.00000391.$$

Also the constant $A$ in the error term of $G(x)$ does not make a significant difference to the results.

| $\alpha$ | $a$ | $b$ |
|---|---|---|
| 0.96 | 10 | 25 |
| 0.94 | 9 | 28 |
| 0.92 | 9 | 32 |
| 0.91 | 8 | 34 |
| 0.90 | 8 | 36 |
| 0.88 | 7 | 41 |
| 0.85 | 7 | 52 |
| 0.63093 | 1 | 11500 |

Table 2.1: This table shows that all possible exceptions to $G(p^{\alpha}) > 0$ occur when $a \leq \omega(p-1) \leq b$.

Now in order to prove Question 2.1 we need a way of checking the cases that are left from this stage of the proof. In the next chapter we introduce the sieve which will be the main tool in the next step in answering Question 2.1.

# 3

# Sieving

After the work of the previous chapter we are left with only a finite number of cases to check in order to prove Question 2.1 for various $\alpha$. Now in order to find a way to answer the question for these cases we need to understand why for these particular cases, the main term of (2.16) does not outweigh the error term.

The intervals in Table 2.1 are all for relatively small $\omega(p-1)$. To prove that $G(x) > 0$ for a particular case of $\omega(p-1)$, we have to prove it for all $p \geq p_0$ where $p_0$ is the lower bound for $p$

$$p_0 = \max\left\{2.5 \times 10^{15},\ 1 + 2 \cdot 3 \cdots 5 \cdots p_{\omega(p-1)}\right\}.$$

When $\omega(p-1)$ is small then $p_0$ is small and since the error term in (2.16) decreases as $p_0$ increases, we have that the main term is unable to outweigh the error term. Now if there are some large primes dividing $p-1$ then $p_0$ will increase which will decrease the error term in (2.16), making it more likely that $G(x) > 0$ for this particular case. Therefore if we can construct a version of $G(x)$ that depends on the primes dividing $p-1$, then we should be able to prove Question 2.1 for more cases of $\omega(p-1)$.

In §3.2 of this chapter we will introduce the prime sieve which will result in a version of (2.16) that depends on the primes dividing $p-1$. This results in the reduction of cases left to prove. The prime sieve is frequently used to obtain results on the least primitive root modulo a prime. Not only was the sieve used in Cohen and Trudgian's paper [6] on this topic, it was also used to prove Grosswald's conjecture on the Generalised Riemann Hypothesis [21]. Theorem 2.1 was also proved using a version of the sieve as well as results on the sums of primitive roots [5] and consecutive primitive roots [4]. The following section §3.3 will outline how we will computationally reduce these cases even further. However first we need some background on $e-$free integers that is essential to defining the sieve.

## 3.1 e-free integers

**Definition 3.1.** Let $p$ be a prime and let $e$ be a divisor of $p-1$. Suppose $p \nmid n$ then $n$ is **e−free** if, for any divisor $d$ of $e$, such that $d > 1$, $n \equiv y^d \pmod{p}$ is insoluble.

Note that $2-$free is not the same as square-free.

**Definition 3.2.** Let $p_i$, where $i = 1, \ldots, \omega(n)$ are the distinct prime divisors of $n$. The **radical** of $n$ is

$$\mathrm{Rad}(n) = p_1 p_2 \cdots p_{\omega(n)}.$$

The following proposition is important to the proof of Lemma 3.1, the indicator function for $e-$free integers. The proposition shows that the definition of $e-$free depends only on the distinct prime divisors of $e$. This proposition leads us to prove Proposition 3.2, which shows how $e-$free integers relate to primitive roots. Not only does Proposition 3.1 allow us to prove this important fact, it allows us to find the link between the prime divisors of $p-1$ and the error term of (2.16). This link is essential to the next step in answering Question 2.1 and will be outlined in §3.2.

**Proposition 3.1.** *$n$ is $e-$free $\iff$ $n$ is $Rad(e)-$free.*

*Proof.* Let $p$ be a prime and let $e \mid p-1$.
Clearly $e-$free $\implies$ $\mathrm{Rad}(e)-$free as all divisors of $\mathrm{Rad}(e)$ also divide $e$ and are not equal to 1.
Now to prove that $\mathrm{Rad}(e)-$free $\implies$ $e-$free assume $n$ (indivisible by $p$) is $\mathrm{Rad}(e)-$free.
That is for all $f \mid \mathrm{Rad}(e)$, $n \equiv x^f \pmod{p}$ is insoluble (where $x$ is arbitrary).
Suppose $n$ is not $e-$free, then there exists $d \mid e$ with $d > 1$ such that

$$n \equiv y^d \pmod{p} \quad \text{for some } y.$$

Since $d \mid e$, there exists $f \mid \text{Rad}(e)$ such that $f \mid d$, in particular there exists an integer $a$ such that $fa = d$. Therefore,

$$n \equiv y^d \equiv y^{fa} \equiv x^f \pmod{p}, \quad \text{where } x = y^a.$$

This is a contradiction to the assumption that $n$ is $\text{Rad}(e)-$free. □

**Proposition 3.2.** *$n$ is a primitive root modulo $p \iff n$ is $(p-1)-$free.*

*Proof.* Fix $g$ a primitive root modulo $p$.
First we prove that if $n$ is $(p-1)-$free then $n$ is a primitive root modulo $p$.
Let $n \equiv g^k \pmod{p}$ then by Remark 1.1 it suffices to show

$$n \text{ is } (p-1) - \text{free} \Longrightarrow (k,\, p-1) = 1.$$

Assume $n$ is $(p-1)-$free, that is $n \not\equiv y^d \pmod{p}$, for all $d \mid p-1$, with $d > 1$.
Suppose $(k,\, p-1) \neq 1$ in particular $(k,\, d_0) = D$ for some $d_0 \mid p-1$. Then $D \mid k$ and $D \mid d_0$, in particular there exists integer $a$ such that $k = aD$. Therefore,

$$n \equiv g^k \equiv g^{aD} \equiv y^D \pmod{p}, \quad \text{where } y = g^a.$$

However $D \mid d_0$ and $d_0 \mid p-1$ so $D \mid p-1$, this is a contradiction to $n$ being $(p-1)-$free. Therefore $(k,\, p-1) = 1$.
Next we prove that if $n$ is a primitive root modulo $p$ then $n$ is $(p-1)-$free.
Assume $n$ is a primitive root modulo $p$ then $n \equiv g^k \pmod{p}$ if and only if $(k,\, p-1) = 1$.
Suppose $n$ is not $(p-1)-$free, then there exists $d \mid p-1$ with $d > 1$ such that $n \equiv y^d \pmod{p}$.
Since $g$ is a primitive root modulo $p$, there exists $a \in \{1, \ldots, p-1\}$ such that $y \equiv g^a \pmod{p}$, therefore

$$n \equiv y^d \equiv g^{ad} \pmod{p}.$$

Since $d \mid (p-1)$, $(ad,\, p-1) \neq 1$ which is a contradiction to $n$ being a primitive root modulo $p$. Therefore $n$ is $(p-1)-$free. □

Now with these two propositions we are able to define an indicator function for $e-$free integers. This function is of a very similar form to (2.1), the indicator function for primitive roots modulo a prime. This is expected as Proposition 3.2 shows an equivalence between primitive roots modulo a prime and $(p-1)-$free integers. As we mentioned at the beginning of the chapter we will be defining a sieving version of (2.16) that depends on the primes dividing $p-1$ in the next section. The indicator function for $e-$free integers is an important in the proof of this.

**Lemma 3.1** (Lemma 2 from [21]). *Let $p$ be a prime and let $e$ be a divisor of $p-1$. Then*

$$\frac{\phi(e)}{e}\sum_{d|e}\frac{\mu(d)}{\phi(d)}\sum_{\chi\in\Gamma_d}\chi(n)=\begin{cases}1 & \text{if } n \text{ is } e\text{-free},\\ 0 & \text{otherwise},\end{cases} \tag{3.1}$$

*where $\Gamma_d$ is the set of Dirichlet characters modulo $p$ of order $d$.*

*Proof.* Fix $g$ a primitive root modulo $p$ and let $\mathrm{ind}_g n = v$.
Now let $\chi(n)\in\Gamma_d$ then as $\chi(n)$ is multiplicative we have

$$\chi(n)=\chi(g)^v=e^{2\pi ikv/d} \quad \text{where } k=1,\,\ldots,\,d \text{ and } (k,\,d)=1.$$

So the innermost sum of 3.1 becomes

$$\sum_{\chi\in\Gamma_d}\chi(n)=\sum_{\substack{1\le k\le d\\(k,\,d)=1}}e^{2\pi ikv/d}. \tag{3.2}$$

From Proposition 1.4 we have

$$\sum_{t|(k,\,d)}\mu(t)=\begin{cases}1 & \text{if } (k,\,d)=1,\\ 0 & \text{otherwise.}\end{cases}$$

Using this property we can rewrite (3.2) as the following

$$\sum_{\substack{1\le k\le d\\(k,\,d)=1}}e^{2\pi ikv/d}=\sum_{\substack{1\le k\le d\\t|(k,\,d)}}\mu(t)e^{2\pi ikv/d}.$$

Now since $t\mid(k,\,d)$ we have $t\mid k$ and $t\mid d$ so there exists an integer $a$ such that $k=at$. Therefore

$$\sum_{\substack{1\le k\le d\\t|(k,\,d)}}\mu(t)e^{2\pi ikv/d}=\sum_{\substack{1\le k\le d\\t|k,\,t|d}}\mu(t)e^{2\pi ikv/d}=\sum_{t|d}\mu(t)\sum_{1\le at\le d}e^{2\pi iatv/d}.$$

Hence

$$\sum_{\chi\in\Gamma_d}\chi(n)=\sum_{t|d}\mu(t)\sum_{a=1}^{d/t}e^{2\pi iatv/d}.$$

Suppose $d\nmid tv$ then $e^{2\pi iavt/d}\neq 1$ for all $a$ and so we can consider the following geometric series

$$\sum_{a=1}^{d/t}X^a \quad \text{where } X=e^{2\pi ivt/d}.$$

We therefore have

$$\sum_{a=1}^{d/t} X^a = \frac{X(1-X^{d/t})}{1-X} = 0 \qquad \text{as } X^{d/t} = \left(e^{2\pi i v t/d}\right)^{d/t} = 1.$$

Now consider when $d \mid tv$ then $e^{2\pi i a v t/d} = 1$ for all $a$ and so

$$\sum_{a=1}^{d/t} X^a = \frac{d}{t}.$$

Therefore

$$\sum_{\chi \in \Gamma_d} \chi(n) = \begin{cases} \sum_{t|d} \mu(t) \frac{d}{t} & \text{if } d \mid tv, \\ 0 & \text{otherwise.} \end{cases}$$

Let $m = (d, v)$ $d \mid tv$ and recall Remark 1.2 which states $\mu(ab) = \mu(a)\mu(b)$ when $(a, b) = 1$, as $\mu$ is multiplicative, and $\mu(ab) = 0$ otherwise. Then

$$\begin{aligned} \sum_{\chi \in \Gamma_d} \chi(n) = \sum_{t|m} \mu\left(t \cdot \frac{d}{m}\right) \frac{m}{t} &= \sum_{\substack{t|m \\ (t, \frac{d}{m})=1}} \mu(t) \mu\left(\frac{d}{m}\right) \frac{m}{t} \\ &= \mu\left(\frac{d}{m}\right) m \sum_{\substack{t|m \\ (t, \frac{d}{m})=1}} \mu(t) \frac{1}{t}. \end{aligned}$$

Using the Euler product (Theorem 285 in [15]) we obtain

$$\begin{aligned} \sum_{\chi \in \Gamma_d} \chi(n) &= \mu\left(\frac{d}{m}\right) m \prod_{\substack{p|m \\ p \nmid \frac{d}{m}}} \left(1 + \frac{\mu(p)}{p} + \frac{\mu(p^2)}{p^2} + \frac{\mu(p^3)}{p^3} + \dots\right) \\ &= \mu\left(\frac{d}{m}\right) m \prod_{\substack{p|m \\ p \nmid \frac{d}{m}}} \left(1 - \frac{1}{p}\right). \end{aligned}$$

Using Proposition 1.5 we obtain

$$\sum_{\chi \in \Gamma_d} \chi(n) = \mu\left(\frac{d}{m}\right) \frac{\phi(d)}{\phi\left(\frac{d}{m}\right)}.$$

Consider the case where $n$ is $e-$free, we will show that $(v, d) = 1$ for all $d \mid e$. Suppose $(v, d_0) = D$ for some $d_0 \mid e$, then $D \mid v$ and $D \mid d_0$. So there exists some integer $a$ such that $v = aD$ and so $n = (g^a)^D$ however $D \mid d_0 \mid e$. This is a contradiction to $n$ being $e-$free. Hence when $n$ is $e-$free, $(v, d) = 1$ for all $d \mid e$. For this case consider the following

$$\sum_{d|e} \frac{\mu(d)}{\phi(d)} \sum_{\chi \in \Gamma_d} \chi(n) = \sum_{d|e} \frac{\mu(d)}{\phi(d)} \left(\mu(d) \frac{\phi(d)}{\phi(d)}\right) = \sum_{d|e} \frac{\mu(d)^2}{\phi(d)}.$$

Since $\phi(d)$ and $\mu(d)$ are both multiplicative, the sum over the divisors is multiplicative in $e$ and we obtain

$$\sum_{d|e} \frac{\mu(d)^2}{\phi(d)} = \left(\sum_{d|p_1^{\alpha_1}} \frac{\mu(d)^2}{\phi(d)}\right) \cdots \left(\sum_{d|p_r^{\alpha_r}} \frac{\mu(d)^2}{\phi(d)}\right)$$

where $e = p_1^{\alpha_1} p_2^{\alpha_2} \dots p_r^{\alpha_r}$ is the prime decomposition of $e$. Since $\mu(p^\alpha) = 0$ for $\alpha > 1$ and $\mu(p^\alpha) = -1$ when $\alpha = 1$

$$\sum_{d|e} \frac{\mu(d)^2}{\phi(d)} = \left(1 + \frac{1}{\phi(p_1)}\right) \cdots \left(1 + \frac{1}{\phi(p_r)}\right) = \prod_{p|e} \left(\frac{p}{p-1}\right).$$

Now by Proposition 1.5 we have

$$\sum_{d|e} \frac{\mu(d)}{\phi(d)} \sum_{\chi \in \Gamma_d} \chi(n) = \frac{e}{\phi(e)}.$$

Next consider the case when $n$ is not $e-$free, then $M = (e, v) > 1$. We may assume $e$ is square-free by 3.1. Then

$$\begin{aligned}
\sum_{d|e} \frac{\mu(d)}{\phi(d)} \sum_{\chi \in \Gamma_d} \chi(n) &= \sum_{d|e} \mu(d) \left( \mu\left(\frac{d}{M}\right) \frac{1}{\phi\left(\frac{d}{M}\right)} \right) \\
&= \sum_{d|\frac{e}{M}} \sum_{k|M} \frac{\mu(dk)\mu(d)}{\phi(d)} \\
&= \sum_{d|\frac{e}{M}} \frac{\mu(d)^2}{\phi(d)} \sum_{k|M} \mu(k) = 0.
\end{aligned}$$

The inner sum is zero for all $M > 1$ (Proposition 1.4). Hence

$$\frac{\phi(e)}{e} \sum_{d|e} \frac{\mu(d)}{\phi(d)} \sum_{\chi \in \Gamma_d} \chi(n) = \begin{cases} 1 & \text{if } n \text{ is } e-\text{free}, \\ 0 & \text{otherwise.} \end{cases}$$

□

Now that we have an indicator function for $e-$free integers, we can sum this over the square-free integers to obtain an expression for the number of square-free and $e-$free integers. Given $x$ and prime $p$ such that $x < p$, let $N_e^\square(p, x)$ denote the number of square-free and $e-$free integers less than $x$. Then from Proposition 3.2 we have $N^\square(p, x) = N_{p-1}^\square(p, x)$. Now summing the indicator function for $e-$free integers (3.1) over the square-free integers less than $x$ we have

$$N_e^\square(p, x) = \frac{\phi(e)}{e} \left\{ \sum_{\substack{n \le x \\ n = \square-\text{free}}} 1 + \sum_{\substack{d|e \\ d>1}} \frac{\mu(d)}{\phi(d)} \sum_{\chi \in \Gamma_d} \sum_{\substack{n \le x \\ n = \square-\text{free}}} \chi(n) \right\}. \tag{3.3}$$

In the next section we will use (3.3) to obtain a sieving inequality. This inequality will enable us to define a sieving version of $G(x)$ (2.16). This will allow us to answer Question 2.1 for a number of the cases remaining from §2.1.

## 3.2 Results with sieving

The next step in answering Question 2.1 involves a sieving inequality, as mentioned above. In §2.1 the main term of (2.16) is unable to outweigh the error term for a fixed $\omega(p-1)$ when $p$ is small. Therefore if there are some large primes dividing $p-1$, $p$ will be large and the main term is more likely to outweigh the error term, and our theorem will be proved for that specific case. The sieve uses this idea, by taking into account the primes dividing $p-1$.

Given prime $p$, let $k$ be a divisor of $\text{Rad}(p-1)$. Then write

$$\text{Rad}(p-1) = kp_1 \cdots p_s \tag{3.4}$$

where $p_1, \ldots, p_s$ are distinct primes with $1 \leq s \leq \omega(p-1)$ and $k$ is the product of the smallest $\omega(p-1) - s$ distinct primes dividing $p-1$. This is called sieving with core $k$ and $s$ sieving primes. Recall $N_e^{\square}(p, x)$ denote the number of integers less than $x$ that are both square-free and $e-$free (given $p$) and $N^{\square}(p, x)$ denotes the number of square-free primitive roots modulo $p$ less than $x$. The following lemma defines an inequality relating the number of square-free primitive roots with the number of square-free and $e-$free integers.

**Lemma 3.2** (Lemma 2 from [6]). *Given a prime $p$, assume that (3.4) holds. Then*

$$N^{\square}(p, x) \geq \sum_{i=1}^{s} N_{kp_i}^{\square}(p, x) - (s-1)N_k^{\square}(p, x) \tag{3.5}$$

$$= \sum_{i=1}^{s} \left\{ N_{kp_i}^{\square}(p, x) - \left(1 - \frac{1}{p_i}\right) N_k^{\square}(p, x) \right\} + \delta N_k^{\square}(p, x) \tag{3.6}$$

*where*

$$\delta = 1 - \sum_{i=1}^{s} \frac{1}{p_i}. \tag{3.7}$$

*Proof.* Given $n$, a square-free and $k-$free integer, we will show that the right hand side of (3.5) contributes 1 if $n$ is additionally $p_i-$free for all $i$, and otherwise contributes a non-positive quantity.

Suppose for all $i$, $n$ is additionally $p_i-$free then it contributes $s$ to $\sum_{i=1}^{s} N_{kp_i}^{\square}(p, x)$

and $(s-1)$ to $(s-1)N_k^{\square}(p,\, x)$. Therefore $n$ contributes $s-(s-1)=1$ to the right hand side of (3.5). Now Suppose $n$ is not additionally $p_i-$free for all $i$ then it contributes $t$ to $\sum_{i=1}^{s} N_{kp_i}^{\square}(p,\, x)$, where $0 \le t < s$ and again $(s-1)$ to $(s-1)N_k^{\square}(p,\, x)$. In this case $n$ contributes $t-(s-1)$ which is negative as $t<s$.
By the definition of $e-$free, if an integer $x$ is $k-$free and $p_i-$free, for all $i$, then $x$ is $\text{Rad}(p-1)-$free. Furthermore by Proposition 3.1 $x$ is $(p-1)-$free and therefore, by Proposition 3.2, $x$ is a primitive root modulo $p$.
Hence for a square-free and $k-$free integer $n$, if $n$ is additionally $p_i-$free then it is a square-free primitive root modulo $p$. Hence we obtain the inequality (3.5).
Now to obtain (3.6) consider

$$\begin{aligned}\sum_{i=1}^{s} N_{kp_i}^{\square}(p,\, x) - \left(1-\frac{1}{p_i}\right) N_k^{\square}(p,\, x) &= \sum_{i=1}^{s} N_{kp_i}^{\square}(p,\, x) - sN_k^{\square}(p,\, x) + N_k^{\square}(p,\, x)\sum_{i=1}^{s}\frac{1}{p_i}\\ &= \sum_{i=1}^{s} N_{kp_i}^{\square}(p,\, x) - sN_k^{\square}(p,\, x) + N_k^{\square}(p,\, x)(1-\delta)\\ &= \sum_{i=1}^{s} N_{kp_i}^{\square}(p,\, x) - (s-1)N_k^{\square}(p,\, x) - \delta N_k^{\square}(p,\, x).\end{aligned}$$

□

To find a lower bound on $N_k^{\square}(p,\, x)$ we use the same method outlined in §2.1. Recall (2.13), the bound on the sum of Dirichlet characters, and (2.14), the lower bound for the number of square-free integers less than $x$. These bounds can be substituted into the indicator function (3.1) along with the number of square-free divisors of $k$

$$\sum_{\substack{d|k\\ d>1}} |\mu(d)| = 2^{\omega(k)} - 1$$

to obtain

$$N_k^{\square} \ge \frac{\phi(k)}{k}\left(x^{1/2}p^{1/4}\right)\left\{\frac{6}{\pi^2}x^{1/2}p^{-1/4} - \frac{0.104}{p^{1/4}} - \left(2^{\omega(k)}-1\right)E\right\}. \tag{3.8}$$

Recall

$$E \le (\log p)^{1/2}\left(\frac{6}{\pi^2}c^{1/2}(1+\hat{D}) + \frac{10}{3}Ap^{1/8-\alpha/4}c^{3/4}(\log p)^{1/4}\right).$$

A similar bound is found for each $N_{kp_i}^{\square}$ where $i=1,\, \ldots,\, s$ and since

$$\frac{\phi(kp_i)}{kp_i} = \frac{\phi(k)}{k}\frac{p_i-1}{p_i} = \left(1-\frac{1}{p_i}\right)\frac{\phi(k)}{k}$$

we have

$$\left|N^{\square}_{kp_i}(p,\,x) - \left(1-\frac{1}{p_i}\right)N^{\square}_{k}(p,\,x)\right|$$
$$\leq \left(1-\frac{1}{p_i}\right)\frac{\phi(k)}{k}(x^{1/2}p^{1/4})\left\{\left(2^{\omega(kp_i)}-1\right)E-\left(2^{\omega(k)}-1\right)E\right\}$$
$$= \left(1-\frac{1}{p_i}\right)\frac{\phi(k)}{k}(x^{1/2}p^{1/4})\left(2^{\omega(kp_i)}-2^{\omega(k)}\right)E.$$

Since $\omega(n)$ is the number of distinct prime factors of $n$

$$2^{\omega(kp_i)}-2^{\omega(p_i)}=2^{\omega(k)+1}-2^{\omega(k)}=2^{\omega(k)}2-2^{\omega(k)}=2^{\omega(k)}.$$

Therefore

$$\left|N^{\square}_{kp_i}(p,\,x)-\left(1-\frac{1}{p_i}\right)N^{\square}_{k}(p,\,x)\right| \leq \left(1-\frac{1}{p_i}\right)\frac{\phi(k)}{k}(x^{1/2}p^{1/4})2^{\omega(k)}E. \tag{3.9}$$

The following theorem defines the sieving version of (2.16).

**Theorem 3.1.** *Given the prime $p$, assume (3.4) holds. Let $\delta$ be given by (3.7). Suppose $\delta > 0$ and let*

$$\Delta = \frac{s-1}{\delta}+2.$$

*If*

$$G_s(x) := x^{1/2}p^{-1/4}-\frac{\pi^2}{6}\left(\frac{0.104}{p^{1/4}}+\left(2^{\omega(k)}\Delta+1\right)E\right)>0 \tag{3.10}$$

*then $p$ has a square-free primitive root less than $x$.*

*Proof.* Note

$$\sum_{i=1}^{s}\left(1-\frac{1}{p_i}\right)=s-\sum_{i=1}^{s}\frac{1}{p_i}=s-1+\delta.$$

Substituting (3.9) and (3.8) into Lemma 3.2 and using the above property we obtain

$$N^{\square}(p,\,x)\geq -\sum_{i=1}^{s}\left(1-\frac{1}{p_i}\right)\frac{\phi(k)}{k}x^{1/2}p^{1/4}2^{\omega(k)}E$$
$$+\,\delta\frac{\phi(k)}{k}x^{1/2}p^{1/4}\left\{\frac{6}{\pi^2}x^{1/2}p^{-1/4}-\frac{0.104}{p^{1/4}}-\left(2^{\omega(k)}-1\right)E\right\}$$
$$=\delta\frac{\phi(k)}{k}x^{1/2}p^{1/4}\left\{2^{\omega(k)}E\left(-\frac{1}{\delta}(s-1+\delta)-1\right)+E+\frac{6}{\pi^2}x^{1/2}p^{-1/4}-\frac{0.104}{p^{1/4}}\right\}$$
$$=\delta\frac{\phi(k)}{k}x^{1/2}p^{1/4}\left\{2^{\omega(k)}E\left(-\frac{s-1}{\delta}-2\right)+E+\frac{6}{\pi^2}x^{1/2}p^{-1/4}-\frac{0.104}{p^{1/4}}\right\}$$
$$=\delta\frac{\phi(k)}{k}x^{1/2}p^{1/4}\left\{\frac{6}{\pi^2}x^{1/2}p^{-1/4}-\frac{0.104}{p^{1/4}}-\left(2^{\omega(k)}\Delta+1\right)E\right\}.$$

Therefore $N^{\square}(p, x) > 0$ if

$$x^{1/2}p^{-1/4} - \frac{\pi^2}{6}\left(\frac{0.104}{p^{1/4}} + \left(2^{\omega(k)}\Delta + 1\right)E\right) > 0.$$

□

Note that $G_s(x)$ is very similar to $G(x)$ from §2.1, however now $G_s(x)$ depends on the primes dividing $p-1$. The sieve was introduced because we needed $G(p^\alpha) > 0$ for primes such that $\omega(p-1)$ is small. As we discussed in the beginning of this section, if there are some large primes dividing $p-1$ then $p_0$ will increase, for fixed $n = \omega(p-1)$, and the error term in (2.16) will decrease. The error term in (3.10) depends on $\delta$ which in turn depends on the primes dividing $p-1$. The primes dividing $p-1$ affect the error term just as we expect: if there are some large primes dividing $p-1$ then $\delta$ will increase, this results in the error term of (3.10) decreasing. This means for that particular case the main term is more likely to out weigh the error term.

Note that for any particular case $\omega(p-1) = n$, we can choose the number of sieving primes, $s$. That is, we can choose $s$, for a fixed $n$, that minimises the error term in (3.10). At first one might think that we should just choose $s = \omega(p-1)$ because then we are including all primes in $\delta$ and $2^{\omega(k)}$ will be minimised. However $\delta > 0$ and also the numerator of $\Delta$ would be maximised. So choosing $s = \omega(p-1)$ will not necessarily minimise the error term. To find the right balance between $k$ and $s$ for a fixed $n$ the optimal value of $s$ is chosen such that it minimises the error term in (3.10).

To prove Question 2.1 with $\alpha = 0.9$ the cases left to check are $8 \le \omega(p-1) \le 36$. Consider the first case, $\omega(p-1) = 8$. We can find the optimal $s$ by calculating the error term from (3.10) for each $1 \le s \le 8$ and selecting the $s$ that generates the smallest error term. For this case $s = 6$ is optimal. As described above we would like $\delta$ to be large and so the worst case is if $k = 2 \cdot 3$ and the sieving primes are the next 6 smallest primes. For this case $\delta \ge 1 - (5^{-1} + 7^{-1} + \cdots + p_8^{-1})$. Now substituting $\delta$, $s = 6$, $x = p^{0.9}$ and $\omega(k) = 2$ into (3.10) shows for $\omega(p-1) = 8$

$$G_s(p^{0.9}) > 0 \text{ for all } p > 1.42 \times 10^{13}.$$

Also from Theorem 2.1 as well as the computation mentioned in §2.1, we have that all primes less than $2.5 \times 10^{15}$ have a square-free primitive root less than $p^{0.9}$. Therefore Question 2.1 with $\alpha = 0.9$ is true for $\omega(p-1) = 8$.

Next consider the last case, $\omega(p-1) = 36$. The optimal $s$ for this case is 32 and so $k = 2 \cdot 3 \cdot 5 \cdot 7$ and $\delta \geq 1 - (11^{-1} + 13^{-1} + \cdots + p_{36}^{-1})$. Substituting $\delta$, $s = 32$ and $\omega(k) = 4$ into (3.10) shows that

$$G_s(p^{0.9}) > 0 \text{ for all } p > 2.98 \times 10^{20}.$$

Since $\omega(p-1) = 36$, $p \geq 1+2\cdot3\cdot5\cdots p_{36} = 2.25\times10^{59}$. Therefore Question 2.1 with $\alpha = 0.9$ is also true for $\omega(p-1) = 36$. Continuing in this way for each remaining case, $9 \leq \omega(p-1) \leq 35$, determining the optimal $s$ and using it to find when (3.10) is true, we can prove Question 2.1 with $\alpha = 0.9$ for all cases except $\omega(p-1) = 13$.

Table 3.1 shows the cases of Question 2.1 that cannot be answered at this stage. Compared to Table 2.1 we can see the improvement that the sieved equation makes. By using Theorem 3.1 we have significantly decreased the number of cases left to answer Question 2.1. For the lower bound on the exponent $\alpha$ we have gone from 11500 cases of $\omega(p-1)$ left to prove at the end of §2.1 to just 39 cases left to prove. Equivalently we have proved the following,

$$g^{\square}(p) < p^{0.63093} \text{ for all } p > 9.63 \times 10^{65}.$$

We have also shown that without the computational algorithm outlined in the next section we can prove the following.

*All primes $p$ have a primitive root less than $p^{0.91}$.*

At the same stage in [6], Cohen and Trudgian were unable to prove Question 2.1 with $\alpha = 0.96$ for $\omega(p-1) = 13$. For this particular case they needed to use the next step outlined in §3.3. This improvement from [6] shows how the changes we made to the bounds in §2.1 have carried through to improve (3.10).

| $\alpha$ | $a$ | $b$ |
|---|---|---|
| 0.96 | – | – |
| 0.94 | – | – |
| 0.92 | – | – |
| 0.91 | – | – |
| 0.90 | 13 | 13 |
| 0.88 | 11 | 14 |
| 0.85 | 9 | 15 |
| 0.63093 | 1 | 39 |

Table 3.1: This table shows that all possible exceptions to $G_s(p^\alpha) > 0$ occur when $a \leq \omega(p-1) \leq b$.

## 3.3 Results using the prime divisor tree

This section outlines the final step in our proof of Question 2.1. As explained above Theorem 3.1 depends on the primes dividing $p-1$. If there are some large primes dividing $p-1$ then $\delta$ will increase and in turn the error term in (3.10) will decrease. Therefore we would like to take advantage of when there are some large primes dividing $p-1$. To do this we divide each case of $\omega(p-1)$ into sub cases depending on the primes dividing $p-1$. We can do this using an algorithm that works through a prime divisor tree.

For example consider the case $\alpha = 0.9$ and $\omega(p-1) = 13$. After the second step in our proof, described in §3.2, we have that all possible exceptions for Question 2.1 to be true occur when $p \in [2.5 \times 10^{15}, 4.17 \times 10^{15}]$.

Since we know that 2 divides $p-1$ this is our base case (or base node) with the above interval of possible exceptions. The base node has $s = 10$ and therefore $\delta \geq 1 - (7^{-1} + 11^{-1} + \cdots + p_{13}^{-1}) = 0.416$.

Now we look at the next level of our tree which has two nodes, when $3 \mid p-1$ and when $3 \nmid p-1$. First consider the case when $3 \nmid p-1$. Choosing $s = 10$ gives $k = 2 \cdot 5 \cdot 7$ and $\delta \geq 1 - (11^{-1} + 13^{-1} + \cdots + p_{14}^{-1}) = 0.536$. We can see that $\delta$ has increased from the base node and so the error term in (3.10) will decrease, in turn decreasing the upper bound on the interval of possible exceptions. Now substituting $s$ and $\delta$ into (3.10) we have

$$G(p^{0.9}) > 0 \text{ for all } p > 1.27 \times 10^{15}.$$

Also $p > \max\left\{2.5 \times 10^{15},\, 1 + 2 \cdot 5 \cdot 7 \ldots p_{14}\right\} = 2.5 \times 10^{15}$. Therefore when $\omega(p-1) = 13$ and $3 \nmid p-1$, Question 2.1 with $\alpha = 0.9$ is true. Therefore we do not have to explore any more of the prime divisor tree on this branch, however we still need to look at the other side of the tree where $3 \mid p-1$. In this case we know that both 2 and 3 divide $p-1$ and therefore 6 divides $p-1$. This means that all possible exceptions are of the form $p = 2 \cdot 3 \cdot m \in [2.5 \times 10^{15},\, 4.17 \times 10^{15}]$. However there are many integers of this form and searching for a square-free primitive root of all primes of this form would be very time consuming. Hence this node branches into two nodes, $5 \mid p-1$ and $5 \nmid p-1$.

Now we continue in the same way as we did for the previous level of the tree. There are three possible outcomes at each node. Outcome one is that the interval of possible exceptions is empty ($3 \nmid p-1$ node from above) and so there is no further branching from the node. Outcome two is when the number of possible exceptions in the interval is too large to search exhaustively ($3 \mid p-1$ node from above) and so the node branches into two nodes. Outcome three is that the number of possible exceptions in the interval is small enough to search exhaustively. In this case we take all the integers of the correct form (with the right divisors and non-divisors for that node) from the interval and eliminate those that are not prime and that do not have $\omega(p-1) = 13$. Then for the remaining list of possible exceptions we verify that each prime $p$ in the list has a square-free primitive root less than $p^{0.9}$.

The tree continues until all possible exceptions have been searched exhaustively and therefore there are no possible exceptions when $\omega(p-1) = 13$. Since this is the only case to check when $\alpha = 0.9$, we are done. If there were more cases to check the above method would be applied to each case of $\omega(p-1)$.

The next section describes the algorithm that implements the above tree for a given exponent $\alpha$. The algorithm is extended from the Sage code used in [21].

### 3.3.1 Running the algorithm

To run the algorithm described above the inputs needed are the exponent, $\alpha$, and $n$, the values of $\omega(p-1)$ that need to be checked. Note that $n$ can be a list if there is more than one value to check. The explicit constant $A$ from Lemma 1.1, is also an input, however as discussed in §2.1 we set $A = 0.679091$, as smaller constants do not make a significant difference to the results.

The algorithm runs on a queue of unexplored nodes. We know the first case, $2 \mid p-1$ is true and so the first node that is added to the queue is $(n, 2, \phi, \infty)$. The information stored in this first node is $n$, $\{2\}$ (the list of primes dividing $p-1$), $\{\phi\}$ (the list of primes that do not divide $p-1$) and the upper bound of infinity. If this node branches then an upper bound is computed for this node is and stored in the two resulting nodes, the upper bound is infinity here as the base node did not branch from another node. While there are nodes in the queue (there are still cases to check for possible exceptions), a node is selected and the Square-Free Tree, Algorithm 3.4, is run on this node. If Algorithm 3.4 results in two nodes (a prime does divide and a prime does not divide), these two nodes are added to the queue and the next node is selected. If the node is explored in Algorithm 3.4 (this is outcome three described in the previous section) then there is no node to add to the queue and the next node is selected. The algorithm will run until there are no nodes left in the queue and hence all cases are explored.

The following algorithms describe all the different parts that come together to prove Question 2.1 using the tree described above. First we have the Optimal S algorithm that, at each node will calculate which value of $s$ will minimise the upper bound. This is an improvement from the algorithm in [21] as we are calculating the optimal value of $s$ at each node (as it will change depending on what primes do not divide $p-1$) instead of using the same $s$ for the whole tree.

**Algorithm 3.1** (Optimal S).

```
def optimalS(n, L):
    # L contains the first n primes not in Y
    best_A = ∞
    p := the product of the first n primes + 1
    for s ∈ [1, n]:
        M := the largest s primes in L
        δ := 1 − (q_1^{-1} + q_2^{-1} + ··· + q_s^{-1}) where q_i ∈ M
        A = the error term of (3.10)
        if δ ≤ 0:
            continue # go to the next iteration
        if A < best_A:
            best_A = A # A is the new lowest error
            best_s = s # the corresponding s to this new low error
    return best_s # returns the optimal s for n and L.
```

In order to find the upper bound, we need to find $p$ such that $G_s(p^\alpha) = 0$ and to do that we use the following two algorithms. The Sage code in [21] used the inbuilt root finder which works for their particular values of $n$ and $\alpha = 0.96$. However as $\alpha$ is decreased, the root finder becomes unstable and will not reliably find a root of $G_s(p^\alpha)$ for $\alpha \leq 0.9$. In order to address this we developed Algorithm 3.2 which while slower is more stable.

As there is some error with the root finder (a precision that can be set) we ensure that the root found is an upper bound (an overestimation) so that the range of possible exceptions is not smaller than it should be. As described in Algorithm 3.3, if the root found is close to the lower bound then the precision is increased and a more precise root is found, this continues until we reach our maximal precision, or the root found is smaller than the lower bound. If the root found is much larger than the lower bound then we do not increase the precision as there is no chance that the bounds will overlap with a higher precision. Another possible case is if the precision is too low and the root is higher than the upper bound of the node from which the current node branched. The upper bound should either stay the same or decrease as we move through the tree.

**Algorithm 3.2** (Zig-Zag Root Finder).

```
def zigzag_root_finding(F, num, inc, inc_div, inc_end):
  while True:
    if inc < 0:
      # Looking for function to become negative
      while F(num) ≥ 0:
        # if the function at num does not reach/pass the root of F, then num is decreased by inc
        # this continues until num passes/undershoots the root of F
        num = num + inc
        if num < 1
          num = max{num, 1}
          # num cannot be less than 1 so the loop is broken and hence the precision is increased and reverses direction.
          break
    else:
      # Looking for function to become positive
      while F(num) ≤ 0:
        # if function at num does not reach/pass the root of F, then num is increases by inc
        # this continues until num passes/overshoots the root of F
        num = num + inc
    if |inc| ≤ inc_end: # if we are within the desired precision
      # check we are overshooting above (F(num) > 0 we have an upper bound on the root)
      if inc > 0:
        break # we have found an upper bound on the root of F within the desired precision
      else:
        # we are currently below the root, flip sign and keep this precision (does one final loop in order to have an upper bound on the root)
        inc = −inc
        continue
    # reverse the direction with smaller incrementation.
    inc = −inc/inc_div
  return num # returns the root of F
```

**Algorithm 3.3** (Find Sign Change).

```
def find_sign_change(n, s, d, lb, lub):
    def F(p) = G_s(p^α)  (3.10) # The function for root finding
    num = lb # lower bound is the starting point
    # Skip root finding if the lower and upper bounds cross
    if F(lb − 1) ≥ 0:
        return lb − 1
    curr_start  = 10^21 # starting precision
    curr_prec  = 10^17 # maximum precision reached
    prec_div = 100
    # find root, starting at num,increment by curr_start, find root
        to within inc_end as described in Algorithm 3.2
    num = zigzag_root_finding(F, num, curr_start, 10, curr_prec)
    curr_start = curr_prec
    curr_prec = curr_prec/prec_div
    # if the lower bound is within the precision of the current
        root, then precision is increased and a more precise root
        is found.
    while num − lb ≤ curr_start and curr_start ≥ max_prec_rf:
        num = zigzag_root_finding(F, num, curr_start, 10,
            curr_prec)
        curr_start = curr_prec
        curr_prec = curr_prec/prec_div
    if lub − num < 0:
        # the root found is larger than the upper bound from the
            parent node
        assert |lub − num| ≤ curr_start # assert this is within
            precision
        return lub # Return the parent ub instead, since it must
            be lower
    return num # returns the upper bound
```

When one of the nodes is selected from the queue as described above, the Algorithm 3.4 is applied to it as described below. This algorithm takes the lists $\{X\}$ and $\{Y\}$, the primes that divide $p-1$ and do not divide $p-1$ respectively, to calculate the optimal $s$ for this particular node. This is then used to determine the new upper bound and lower bound.

Now there are four possible cases that the algorithm deals with. Case 1 is when the lower and the upper bound overlap, therefore there are no exceptions to check and Algorithm 3.4 ends. The second case is when the range of possible exceptions is too large to search exhaustively. Here the tree branches and the two new nodes are added to the queue. The third case is when the range of possible exceptions is small enough to check exhaustively. All the numbers in the range that are divisible by all the primes in $X$ are kept as possible exceptions if they are prime and $\omega(p-1) = n$. Then for these possible exceptions a square-free primitive root is found that is less than $p^\alpha$ using Algorithm 3.5. The last case was introduced while trialling the algorithm, as the tree was unbounded leading to extremely long run times. This was due to some branches that were growing extremely long in the $Y$ direction as each addition of a prime to the list $Y$ did not make a significant decrease to the upper bound if there were already a significant number of primes in $X$.

For example consider $n = 13$ and $s = 10$. If $X = \{2, 3, 5, 7, 13, 17, 23, 29\}$ and $Y = \{11, 19, 29\}$ then $\delta \geq 0.515$. Now if we look at the next branch where $31 \nmid p-1$ we have $\delta \geq 0.529$. We can see here that the addition of the larger prime, 31 to $Y$ does not significantly increase $\delta$, so the upper bound will not be significantly decreased and the node will branch. However if 31 did not decrease the upper bound significantly then any larger primes will not either. This results in the branch growing in the $Y$ direction as the range will not get small enough to enumerate. Note that the lower bound also increases with every prime that is added to $Y$ however the lower bound increases very slowly. At this stage the lower bound increases much slower than we can check the remaining possible exceptions.

To overcome this problem an exhaustive search is forced when the number of primes in $Y$ is twice the number of primes in $X$ and the number of primes in $X$ is most of $n$. Therefore in this case where an exhaustive search is forced on what can be a large interval, the time taken for enumeration can increase, however it allows the algorithm to complete.

When the number of possible exceptions in the interval, that is integers of the form $p = m\prod_{q\in X} q - 1$ in the interval, is small, for example there are less than a million, the node is explored. The node is also explored in the above case when the number of primes not dividing $p-1$ is large and enumeration is forced. When this is the case, first we take all the integers of the above form and use Algorithm 3.6 to remove all those that are divisible by one or more of the primes in $Y$. Then we remove all those that are not prime and do not have $\omega(p-1) = n$. The remaining list of possible exceptions, $p$, is much smaller and so we use Algorithm 3.5 to find a square-free primitive root less than $p^\alpha$ for each $p$.

**Algorithm 3.4** (Square-Free Tree).

```
def square_free_tree(n, X, Y, ub):
    L := the first n primes not in Y
    lower := max{1 + q_1 q_2 ... q_n, 2.5 · 10^15} where q_i ∈ L
    upper := 0
    s = optimalS(n, L) (3.1)
    M := the largest s primes in L
    δ := 1 − (q_1 + q_2 + ··· + q_s) where q_i ∈ M
    p := find_sign_change(n, s, δ, lower, ub)
    upper := p
    assert upper > 0
    prodX := q_1 q_2 ... where q_i ∈ X # product of the primes q | p − 1
    enum := (upper − lower)/prodX # number of possible exceptions in the interval
    if enum ≤ 0:
        there are no exceptions at this node
    elif enum > range_limit and len(X)< ⌈0.8n⌉ and len(Y) < 2len(X):
    # The number of possible exceptions is too large to search exhaustively and enumeration is not forced
        q := the smallest prime not in X ∪ Y
        X1 := X ∩ {q}
        (n, X1, Y, upper) added to the list of unexplored nodes
        Y1 := Y ∩ {q}
        (n, X, Y1, upper) added to the list of unexplored nodes
    else: # the number of possible exceptions is small enough to search exhaustively or enumeration is forced
        I := [(lower − 1)/prodX, (upper − 1)/prodX]
        Using Algorithm 3.6 Remove the integers in I divisible by primes in Y
        for k in I:
          p = k · prodX + 1
          if p is prime and ω(p − 1) = n:
            # now we need to find a square free primitive root of this possible exception
            g := least primitive root of p
            if g is not square-free and g > p^α:
                # use Algorithm 3.5 to find a square-free primitive root of p less than p^a
                if not sfpr(p): # there is no square-free p.r. less than p^a by Algorithm 3.5
                    we have found a counter example # this does not happen
                else there are no exceptions at this node
```

The following algorithm, given a prime $p$, finds a square-free primitive root less than $p^\alpha$. As described above, after filtering out integers that are not prime and do not have $\omega(p-1) = n$ we are left with a small list of possible exceptions. For each of these possible exceptions we first use a built-in Sage function, primitive_root($p$), that finds the least primitive root of $p$. If the primitive root is square-free and less than $p^\alpha$ then Question 2.1 is true for this $p$ and it is removed from the list of possible exceptions. Otherwise, Algorithm 3.5 is applied to the possible exception $p$. The algorithm uses Remark 1.1 which states that if we know one primitive root we can generate all other primitive roots. We start with $g$, the least primitive root of $p$, generated using primitive_root($p$). Then for $k \in [2, p-1]$ such that $(k, p-1) = 1$, $g^k \pmod p$ generates another primitive root. If this primitive root is not less than $p^\alpha$ then the next $k$ is used until it is. Then we test if the primitive root is square-free using is_squarefree($g$), a built-in Sage function, if it is not then we again take the next $k$ until it is. This sounds like a long and expensive way of finding a square-free primitive root less than $p^\alpha$ however usually the least primitive root is square-free and so Algorithm 3.5 does not usually have to be called and if it does, not many cycles of $k$ are used.

**Algorithm 3.5** (Square-Free Primitive Root $< p^\alpha$).

```
def sfpr(p)
    g = primitive_root(p) # find the least square-free p.r. using
        inbuilt function
    found_sqfr = False
    for k ∈ [2, p − 1]:
    # cycle through primitive roots until the p.r. that is found
        is square-free and less than p^α
     if gcd(k, p − 1) = 1:
         g = g^k (mod p)
         if g > p^α:
             continue
         if g is square-free:
             found_sqfr = True
             break
    return found_sqfr # returns True if p has a square-free p.r.
        less than p^α
```

The following algorithm removes all integers in the interval $[a, b]$ that are divisible by at least one prime in $Y$. This is an improvement on the Sieve function used in the Sage code used in [21] as it is less memory intensive.

**Algorithm 3.6** (Sieve).

```
def sieve(a, b, Y)
    for n ∈ [a, b]:
        # check that n is not divisible by any primes in Y.
        div_flag = False #start with no divisors
        for y ∈ Y:
            if n ≡ 0 (mod y)
            div_flag = True # if y | n , break out of loop, if not then move to the next element in Y
            break
        if not div_flag:
            # if n is not divisible by any elements in Y then add n to the list.
            yield n
```

The full Sage code can be found in Appendix A. The algorithm ran in parallel on a Dell Power Edge R110 which has an Intel Xeon E3-1230 v2 with 4 3.2GHz cores and 8GB RAM.

### 3.3.2 Results

The following table shows the output of the algorithm outlined in §3.3.1. After the last stage of the proof of Question 2.1 we have proved the following

*All primes $p$ have a square-free primitive root less than $p^{0.88}$.*

It can be seen that for lower exponent $\alpha$ the prime divisor tree grows very large. For just a 0.02 decrease in $\alpha$ the algorithm takes 60 times as long to complete. This suggests that the computational part of the proof of Question 2.1 does not contribute much improvement to the bound. We also ran the algorithm with $\alpha = 0.85$ however it did not complete in time for this thesis despite running for two weeks on a cluster. It is possible that with more computational power, or just simply more time, the exponent could be decreased.

| $\alpha$ | $n$ | nodes | times checked for p.r. | time |
|---|---|---|---|---|
| 0.90 | 13 | 31 | 964 | 1m14s |
| 0.88 | 11 | 37 | $1.70 \times 10^5$ | 1h03m35s |
| | 12 | 785 | $8.67 \times 10^5$ | |
| | 13 | 683 | $5.74 \times 10^4$ | |
| | 14 | 63 | 244 | |

Table 3.2: $\mathbf{n} = \omega(p-1)$, **nodes** is the number of nodes in the tree, **times checked for p.r.** is the number times we had to find a primitive root less than $p^\alpha$ and **time** is the computation time

# 4

# Conclusion and future work

## 4.1 Conclusion

Using stronger explicit bounds for the number of square-free primitive roots less than $x$, see §1.2.2, and for the Pólya–Vinogradov inequality, see §1.2.4, we are able to improve the best known bound for the least square-free primitive root modulo a prime.

Using the proof outlined in §2.1, §3.2 and §3.3 we have proved the following theorem.

**Theorem 0.1.**

$$g^{\square}(p) < p^{0.88} \quad \textit{for all primes } p.$$

Using stronger bounds on the number of square-free integers and the sum of Dirichlet characters we were able to prove in §2.1 that $g^{\square}(p) < p^{0.96}$ for all primes $p$ with $\omega(p-1) \leq 9$ or $\omega(p-1) \geq 26$. This is an improvement on [6] because without the use of the sieve Cohen and Trudgian were left with cases $8 \leq \omega(p-1) \leq 29$.

In §3.2 we proved, without any computations, that $g^{\square}(p) < p^{0.91}$ for all primes $p$. This result shows that most of the improvement to the bound was from the non-computational methods outlined in this section, while only a slight improvement came from the computations in §3.3. We were also able to prove the following without any computations

**Theorem 0.2.**

$$g^{\square}(p) < p^{0.63093} \quad \textit{for all primes } p < 2.5 \times 10^{15} \textit{ and } p > 9.63 \times 10^{65}.$$

Recall that 0.63093 is the theoretical lower bound for the exponent. In order to prove Theorem 0.2 for all primes, the cases left to check are $1 \leq \omega(p-1) \leq 39$. We are not yet able to prove the result for these cases, however Theorem 0.2 is almost certainly true for all primes.

## 4.2 Future work

As we mentioned above, the bound on the least square-free primitive root modulo a prime still seems to have some room for improvement. One way of lowering the exponent $\alpha$ would be to look for stricter bounds on the sums that are used in the formulation of $G(x)$ (2.16) and $G_s(x)$ (3.10). In particular one could use a stricter bound on the number of square-free integers less than $x$. We came to the conclusion that the slightly improved bounds from Cipu (Lemma 1.1) do not make a significant difference to our results.

Although smaller constants $c$ from the Pólya–Vinogradov inequality (1.3) do exist they often require the Dirichlet characters to have certain properties. However it is unclear that slight improvements to this constant will significantly improve results especially for small $\alpha$. In fact we can actually use Theorem 1.3 as all non-principal Dirichlet characters modulo a prime are primitive (Theorem 8.14 in [1]). This theorem states that for $\chi$ a primitive Dirichlet character modulo $q > 1200$ we have,

$$\left|\sum_{n \leq x} \chi(n)\right| < c\sqrt{q}\log q$$

where

$$c = \left(\frac{2}{\pi^2} + \frac{1}{\log q}\right). \tag{4.1}$$

Using this $c$ as our proof instead of the constant from Theorem 1.2 we are able to prove without any computation that $g^{\square}(p) < p^{0.9}$ for all primes $p$. However for smaller values of $\alpha$, this improved $c$ does not improve the intervals of cases left to check after §3.2 (Table 3.1). Therefore we expect that there will not be a noticeable difference in the computational stage of the proof.

Another way of improving the bound and proving Theorem 0.2, would be to make the code for the computational part of the proof, §3.3, more efficient and run it on a more powerful computer. However our results indicate that only small improvements were made to the bound from computations when compared to the improvements made theoretically.

Perhaps we should rather think back to the reason that square-free primitive roots are of interest. We are interested in square-free primitive roots as bounds on the least prime primitive root are very difficult to lower. Here square-free primitive roots provided a generalisation. Another problem that would be interesting to tackle is the least odd square-free primitive root modulo a prime. All primes are odd square-free integers, except for 2, and we could see this as one step closer to the problem of bounding the least prime primitive root modulo a prime. Recall that the number of square-free integers less than $x$ is asymptotic to $6\pi^{-2}x$. It has also been shown that the number of odd square-free integers less than $x$ is asymptotic to $4\pi^{-2}x$ [17]. This means that in the case of odd square-free primitive roots the main term $G(x)$ (2.16) and $G_s(x)$ (3.10) will be smaller. Since we have proven Theorem 0.2 for a significant number of cases, the problem of odd square-free primitive roots is promising.

Another interesting problem would be to look at the distribution of consecutive square-free primitive roots modulo a prime. This does not lead to more results on the least prime primitive root however it is an interesting problem in itself. This problem is explored by Liu and Dong in [18].

# A

# Sage code

```
import Queue
import multiprocessing
import os
import time
import gc
import numpy as np

gross_list = []

# quad precision real field
def pRR(x):
  return RealField(113)(x)

# Arg params
myA=pRR(0.679091)
mya=pRR(0.88)
my_ns = [11, 12, 13, 14]
trace_file = 'a88.csv'

```

```
# Alg params
write_trace = True
range_leaf_lim = 1e6
force_enum_limit = 5e8
adaptive_root_finding = True
start_pref_rf = pRR(1e21)
default_prec_rf = pRR(1e17)  # Also affects fallback onepass
max_prec_rf = pRR(1e9)
prec_div_rf = pRR(10)
cpus = multiprocessing.cpu_count()
parallel_flag = True
parallel_restart_after = 2500
tol0 = pRR(1e-2)

# so everything prints nicely
def print_trace(msg, n=-1, exp=-1, unexp=-1, exc=-1):
  n_s = '-'
  exp_s = '-'
  unexp_s = '-'
  exc_s = '-'
  pid_s = '%i'%os.getpid()
  if n > 0:
    n_s = '%i'%n
  if exp > 0:
    exp_s = '%i'%exp
  if unexp > 0:
    unexp_s = '%i'%unexp
  if exc > 0:
    exc_s = '%.2g'%exc
  tr_ln = ' p:{: <5.5} | n:{: <5.5} | e:{: <7.7} | u:{: <6.6} | c:{: <7.7} | {:}'
    .format(pid_s, n_s, exp_s, unexp_s, exc_s, msg)
  print(tr_ln)
  if write_trace:
    with open(trace_file, 'a') as fl:
      fl.write('%i,%i,%i,%i,%i,\"%s\"\n' % (os.getpid(), n, exp, unexp, exc, msg))

def sieve(a,b,Y):
  for n in xrange(a,b):
    # Check that n is not divisible by anything in Y
    div_flag = False
    for y in Y:
      if n % y == 0:
        div_flag = True
        break
```

```
    if not div_flag:
      yield n

def pr(n):
  return pRR(prod(prime_range(2,nth_prime(n+1)))+1)

def p0(n):
  return pRR(max(2.5e15, pr(n)))

def c(n):
  return pRR(1/(pi*sqrt(2))+6/(pi*sqrt(2)*log(p0(n)))+1/log(p0(n)))

def d(n,p):
  return pRR(sqrt(p^(mya)/(c(n)*sqrt(p)*log(p))))

def err(n,p):
  err1 = pRR(sqrt(c(n)*log(p))*(1+d(n,p)/(d(n,p)-1))-p^(-0.25))
  err2 = pRR((p^(mya*(-0.25)))*((c(n)*log(p))^0.75)*(p^(1/8))-
    p^(mya*(-0.25))*p^(-.25)/3+(7/3)*(p^(mya*(0.5))*(p^(-0.25))*((d(n,p)-1)^(-3/2))))
  return pRR((6/(pi^2))*err1+myA*err2)

def zigzag_root_finding(F, num, inc, inc_div, inc_end):
  num = max(num, 0)
  while True:
    # Do increment Loop, so we overshoot
    if inc < 0:
      # negative increment, looking for negative sign
      while F(num) >= -tol0:
        num+=inc
        if num < tol0 + 1:
          num = max(num, tol0 + 1)
          # If num < 0, we cannot have this happen, break out so prec is
            # increased/flipped up
          break
    else:
      # positive increment, looking for positive sign
      while F(num) <= tol0:
        num+=inc
    # check if reached desired prec
    if abs(inc) <= inc_end:
      # check we are overshooting above (ie it is currently about to go neg)
      if inc > 0:
        break
      else:
```

```
        # we are currently below zero, flip sign and keep this prec
        inc *= -1
        continue
      # flip sign and reduce increment
      inc /= -inc_div
  return pRR(num)

def find_sign_change(n,s,d,lb,lub,prodX):
  def F(p):
    A = pRR(((pi^2)/6)*(0.104/(p^(0.25))+((2^(n-s))*((s-1)/d+2)+1)*err(n,p)));
    B = pRR((p^(mya/2-0.25)));
    return pRR(B-A);

  if adaptive_root_finding:
    num = pRR(lb)
    # Check now if we need to bother with any of this
    if F(num-1) >= tol0:
      return num-1

    # Start at 1e15 and increase perc by 100 each time
    curr_start = start_pref_rf
    curr_prec = default_prec_rf
    prec_div = prec_div_rf
    # No dowhiles in python so do the first zigzag
    num = zigzag_root_finding(F, num, curr_start, pRR(10), curr_prec)
    curr_start = curr_prec
    curr_prec /= prec_div
    assert F(num) >= tol0
    # Loop to adaptively increase prec if best possible ub at this prec
    # (num - curr_start)
    # gives an enum less than the force_enum_limit
    # Only do this up to a maximum end precision of max_prec_rf
    while (num - curr_start - lb) / prodX <= force_enum_limit and
        curr_prec >= max_prec_rf:
      num = zigzag_root_finding(F, num, curr_start, pRR(10), curr_prec)
      curr_start = curr_prec
      curr_prec /= prec_div
      assert F(num) >= tol0
      # If we are already closer than range_leaf_lim stop increasing prec_div
      if (num - lb) / prodX <= range_leaf_lim:
        break
    # sanity
    assert F(num) >= tol0
    if lub - num < 0:
```

```
      # Maybe it is far from LB so precision didn't go far enough
      assert abs(lub-num) <= curr_start
      # Assert it is within precision of last ubber bound
      # Return the last ub instead, since it must be lower
      return lub
    return num
  else:
    # Check if we even need to bother finding the root
    if F(2.5e15)>1e-3:
      return 2.5e15
    try:
      num = ceil(find_root(F,10,10^100,1, maxiter=75));
      # Increment it up a bit to account for precision
      inc = floor(num/1000)
      while (F(num)<=1):
        num=num+inc
      return num
    except:
      # if it fails just increment through until sign change
      num = 2.5e15
      while F(num)<=1e-2:
        num+=default_prec_rf
      return num

def optimalS(n, Primes=[]):
  if len(Primes) == 0:
    # Primes is default empty list, so make it default first primes
    Primes = prime_range(2,nth_prime(n + 1))
  p = pr(n)
  best_A = infinity
  for s in range(1,n + 1):
    M = Primes[n-s:n] # M is the largest s primes in Primes
    d = 1-RR(sum([1./l for l in M]))
    A = RR(((pi^2)/6)*(0.104/(p^(0.25))+((2^(n-s))*((s-1)/d+2)+1)*err(n,p)))
    if d <= 0:
      continue
    if A < best_A:
      best_A = A
      best_s = s
  return best_s

# Checks for square free pr less than p to mya
def sfpr(p):
  g=primitive_root(p)
```

```
  found_sqfr = False
  for k in xrange(2, p):
    if gcd(k, p-1) == 1:
      g = g^k % p
      if g > p^mya:
        continue
      if is_squarefree(g):
        found_sqfr = True
        break
  return found_sqfr

def SFTree(n,X,Y,ub):
  # We never want X to be longer than n
  assert len(X) <= n
  nodes = list()
  count = 0
  # L contains the first n primes not in Y
  print_trace('  Exploring node (depth: %i)' % (len(X)+len(Y)), n=n)
  L=list(set(prime_range(2,nth_prime(len(Y)+n)+1))-set(Y))
  L.sort()
  prodX=prod(X)
  lower = max(prod(L)+1,2.5*10^15)
  upper = 0
  # Get optimalS at each node now
  s = optimalS(n, L)
  M = L[n-s:n] # M is the largest s primes in L
  d = pRR(1-pRR(sum([pRR(1/pRR(l)) for l in M])))
  p=find_sign_change(n,s,d, lower, ub, prodX)
  upper = p; s1 = s; d1 = d; M1=M
  assert upper > 0
  enum = (upper-lower)/prodX
  print_trace('    Node range: (%.2g, %.2g)' % (lower,upper), n=n)
  print_trace('    Enum: %.2g' % enum, n=n)
  if enum <= 0:
    print_trace("    Empty range. No Children.", n=n)
  elif enum > range_leaf_lim and not (len(X) >= ceil(0.8*n) and
    len(Y) >= 2 * len(X) and enum <= force_enum_limit):
    # q is the smallest prime not in X or Y
    q = next_prime(max(X+Y))
    print_trace("    Getting children with next prime %i"%q, n=n)
    if len(X) < n:
        # Only add another prime divisor if we have less than n prime divisors
        X1=list(X)
        X1.append(q)
```

```
        nodes.append((n,X1,Y,upper))
    Y1=list(Y)
    Y1.append(q)
    nodes.append((n,X,Y1,upper))
  else:
    if enum <= range_leaf_lim:
        print_trace("    Less than %.1g elements in range. Enumerating."
            % range_leaf_lim, n=n)
    else:
        print_trace("    Enumeration forced. %.4g elemens in range." % enum, n=n)
    count = 0
    count_s = 0
    for k in sieve(ceil((lower-1)/prodX),1+floor((upper-1)/prodX),Y):
      p=k*prodX+1;
      if is_prime(p) and len(prime_divisors(p-1))==n:
        count = count + 1;
        # now we need to find a square free primitive root of this exception
        # find a primitive root, if its not square-free then find a square-free one
        g=primitive_root(p)
        if not (is_squarefree(g) and g < p^mya):
          count_s += 1
          if not sfpr(p):
            print_trace("FAIL %i" % p)
            assert 1==0
    print_trace("      Checked " + str(count) + " exceptions.", n=n)
  # Return dict with exceptions and discovered nodes
  gc.collect()
  return {'nodes': nodes, 'exceptions': count}

def para_q(q, info, idle):
  # Loop until queue is empty
  idle[os.getpid()] = False
  while True:
    if info['this_run_nodes'] >= parallel_restart_after:
      idle[os.getpid()] = True
      print_trace('Worker quitting for restart')
      return
    elif q.empty():
      idle[os.getpid()] = True
      if all(idle.values()):
        # wait 1 second and check again
        sleep(1)
        if all(idle.values()):
          # wait one more second and check again
```

```
          # (so all workers are idle twice in 2sec)
          sleep(1)
          if all(idle.values()):
            print_trace('Worker quitting')
            return
    else:
      idle[os.getpid()] = False
      try:
        ret = SFTree(*q.get_nowait())
      except Queue.Empty:
        idle[os.getpid()] = True
        continue
      info['explored_nodes'] += 1
      info['this_run_nodes'] += 1
      for node in ret['nodes']:
        q.put(node)
      info['exceptions'] += ret['exceptions']
      print_trace('  Node explored', n=info['n'],
        exp=info['explored_nodes'], unexp = q.qsize(), exc=info['exceptions'])
      gc.collect()

def para_q_wrap(argtup):
  return para_q(*argtup)

def para_n(n):
  print_trace('STARTING n=%i' %n, n=n)
  manager = multiprocessing.Manager()

  # Set up managed objects
  mq = manager.Queue()
  minfo = manager.dict()
  midle = manager.dict()

  # initialise queue and info
  mq.put((n, [2], [], infinity))
  minfo['explored_nodes'] = 0
  minfo['exceptions'] = 0
  minfo['n'] = n

  arg_l = []
  for _ in range(cpus):
    arg_l.append((mq, minfo, midle))

  while not mq.empty():
```

```
      print_trace('(Re)starting workers', n=n)
      minfo['this_run_nodes'] = 0
      midle.clear()
      pl = multiprocessing.Pool(processes=cpus)
      pl.map(para_q_wrap, arg_l, chunksize=1)
      pl.terminate()
      gc.collect()
  return {'explored': minfo['explored_nodes'],
          'exceptions': minfo['exceptions']}


# Runs in parallel
# Runs a queue of SFTree args, handling stuff returned and adding to queue
def run_q(q, n):
  # Loop until queue is empty
  explored_nodes = 0
  c_exceptions = 0
  while not q.empty():
    ret = SFTree(*q.get())
    explored_nodes += 1
    for node in ret['nodes']:
      q.put(node)
    c_exceptions += ret['exceptions']
    print_trace('  Node explored', n=n, exp=explored_nodes,
        unexp=q.qsize(), exc=c_exceptions)
  return {'explored': explored_nodes,
          'exceptions': c_exceptions}

# Runs not in parallel
def run_n(n):
  print_trace('STARTING n=%i' %n, n=n)
  q = Queue.Queue()
  # Put root on queue
  q.put((n, [2], [], infinity))
  # Call run Queue
  return run_q(q, n)

def run_computation():
  if write_trace:
    with open(trace_file, 'w+') as fl:
      fl.write('pid,n,explored,unexplored,exceptions,message\n')
  start_time = time.time()
  tree_sizes = []
  exceptions_checked = []
```

```
  for n in my_ns:
    if parallel_flag:
      ret = para_n(n)
    else:
      ret = run_n(n)
    tree_sizes.append(ret['explored'])
    exceptions_checked.append(ret['exceptions'])
  print_trace("COMPLETE")
  print_trace("Done for a=%.4f" % mya)
  print_trace("Checked for n in %s" % my_ns)
  print_trace(("Tree sizes: ["+', '.join(['%.4g']*len(tree_sizes))+"]")
    % tuple(tree_sizes))
  print_trace(("Checked exceptions: ["+', '.join(['%.4g']*len(exceptions_checked))+"]")
    % tuple(exceptions_checked))
  s = ceil(time.time() - start_time)
  m = floor(s/60)
  s = s % 60
  h = floor(m/60)
  m = m % 60
  print_trace("Time taken: %dh%02dm%02ds" % (h, m, s))

# DO IT!
if __name__ == '__main__':
  print '--------------------------------------------\n    STARTING COMPUTATION'
  print '--------------------------------------------'
  run_computation()
```